\documentclass[preprint,review,12pt]{elsarticle}

\usepackage{upgreek}

\newcommand{\ADP}{\texttt{ADP}}
\newcommand{\NeurADP}{\texttt{NeurADP}}
\newcommand{\Myopic}{\texttt{Myopic}}
\newcommand{\Time}{\mathcal{T}}
\newcommand{\GraphNetwork}{\mathcal{N}}
\newcommand{\VehiclesSet}{\mathcal{A}}
\newcommand{\RequestsSet}{\mathcal{B}}
\newcommand{\FeasibilitySet}{\mathcal{F}}
\newcommand{\RebalancingFeasibilitySet}{\mathcal{Q}}
\newcommand{\MovementSet}{\mathcal{D}}
\newcommand{\RebalancingPointsSet}{\mathcal{P}}
\newcommand{\DecisionSet}{\mathbf{X}}
\newcommand{\AggregationSet}{\mathcal{G}}
\newcommand{\TimeIntervals}{\texttt{$\delta$}}
\newcommand{\Locations}{\texttt{$\mathcal{L}$}}
\newcommand{\Edges}{\texttt{$\mathcal{E}$}}
\newcommand{\CurrentTime}{\texttt{$t$}}
\newcommand{\StateNote}{\texttt{$S$}}

\newcommand{\VehicleLocation}{\texttt{$loc$}}
\newcommand{\VehicleRequests}{\texttt{$reqs$}}
\newcommand{\IndividualVehicle}{\texttt{$a$}}
\newcommand{\IndividualRequest}{\texttt{$b$}}
\newcommand{\RequestOrigins}{\texttt{$orig$}}
\newcommand{\RequestDestinations}{\texttt{$dest$}}
\newcommand{\RequestCapacity}{\texttt{$pass$}}
\newcommand{\RequestDeadline}{\texttt{$dead$}}
\newcommand{\TravelTime}{\texttt{time}}
\newcommand{\MaxWait}{\texttt{wait\textsubscript{\texttt{max}}}}
\newcommand{\MaxDelay}{\texttt{delay\textsubscript{\texttt{max}}}}
\newcommand{\MaxVehicleCapacity}{\texttt{capacity\textsubscript{\texttt{max}}}}
\newcommand{\MaxVehicleGroups}{\texttt{groups\textsubscript{\texttt{max}}}}
\newcommand{\RebalancingPoint}{p}
\newcommand{\MovementPoint}{d}

\newcommand{\ExogenousInformation}{W}
\newcommand{\NullAction}{\varnothing}
\newcommand{\DecisionTuple}{\textbf{x}}
\newcommand{\PostDecisionState}{S^{\texttt{Vehicle-Post}}}
\newcommand{\FirstTransition}{{\texttt{statepost}}}
\newcommand{\SecondTransition}{{\texttt{statenext}}}
\newcommand{\IndividualVehiclePickedUp}{\texttt{$up$}}
\newcommand{\IndividualVehicleDropOff}{\texttt{$drop$}}
\newcommand{\Reward}{R}
\newcommand{\PostDecision}{\texttt{Post}}

\newcommand{\VehiclePost}{\texttt{Vehicle-Post}}
\newcommand{\EstimatedV}{\bar{v}}
\newcommand{\AggregateLevel}{g}
\newcommand{\AggregationFunction}{M}
\newcommand{\Location}{\ell}
\newcommand{\StepSize}{\alpha}
\newcommand{\Days}{N}
\newcommand{\SingleDay}{n}
\newcommand{\MarginalValue}{\hat{v}}
\newcommand{\Vehicle}{\texttt{Vehicle}}
\newcommand{\Request}{\texttt{Request}}
\newcommand{\OptimalityFunction}{F}

\newcommand{\Weights}{w}
\newcommand{\MatchingVariable}{x}

\usepackage[utf8]{inputenc}
\usepackage{fullpage,times}
\usepackage{amsmath,amssymb}
\usepackage{graphicx}
\usepackage{algorithm}
\usepackage{algorithmicx}
\usepackage{subfig}
\usepackage{caption}
\usepackage{url}
\usepackage[numbers]{natbib} 
\usepackage{placeins}
\usepackage{siunitx}
\usepackage{url}
\usepackage{xspace}
\usepackage[table,xcdraw]{xcolor}

\usepackage{hyperref}
\hypersetup{
colorlinks={true},
linkcolor={Blue3},
urlcolor={Blue3},
citecolor={Blue3},
pdfstartview={FitH},
pdfauthor={Author},
pdftitle={Title},
pdfsubject={Subject}
plainpages=false
}

\usepackage{enumitem}
\DeclareMathOperator*{\argmax}{arg\,max}
\usepackage{booktabs} 
\usepackage{multirow}
\usepackage{multicol}
\usepackage{makecell}

\usepackage[toc,page]{appendix}

\usepackage{setspace}

\let\Algorithm\algorithm

\definecolor{wisconsin-red}{rgb}{0.6,0,0}
\definecolor{darkgreen}{rgb}{0.2,0.6,0.2}
\definecolor{maroon}{rgb}{0.5, 0.0, 0.0}
\definecolor{violet}{rgb}{0.75, 0.0, 1.0}
\definecolor{lightgray}{gray}{0.9}
\definecolor{navyblue}{rgb}{0.0, 0.0, 0.5}
\definecolor{darkmidnightblue}{rgb}{0.0, 0.2, 0.4}
\definecolor{midnightblue}{rgb}{0.0,0.4,0.85}
\definecolor{Gray}{gray}{0.75}
\definecolor{darkgreen}{rgb}{0,0.5,0}
\definecolor{apricot}{rgb}{0.98, 0.81, 0.69}

\newcolumntype{C}[1]{>{\centering\arraybackslash}p{#1}}
\newcolumntype{P}[1]{>{\raggedright\arraybackslash}p{#1}}
\newcolumntype{L}[1]{>{\raggedleft\arraybackslash}p{#1}}

\journal{Elsevier}

\begin{document}\sloppy
\setlength{\parindent}{2em}

\begin{frontmatter}
\title{An Enhanced Approximate Dynamic Programming Approach to On-demand Ride Pooling} 

\author[1]{Arash Dehghan}
\ead{arash.dehghan@torontomu.ca}

\author[1]{Mucahit Cevik\corref{cor1}%
\fnref{fn1}}
\ead{mcevik@torontomu.ca}

\author[2]{Merve Bodur}
\ead{bodur@mie.utoronto.ca}

\cortext[cor1]{Corresponding author}
\fntext[fn1]{Toronto Metropolitan University, Toronto, ON, Canada}
\address[1]{Toronto Metropolitan University, Toronto, ON, Canada}
\address[2]{University of Toronto, Toronto, ON, Canada}

\begin{abstract}
Ride-pooling services have been growing in popularity, increasing the need for efficient and effective operations. The main goal of ride-pooling services is to maximize the number of passengers served while minimizing wait and delay times. However, factors such as the timing and volume of passenger requests, pick-up and drop-off locations, available vehicle capacity, and the trajectory to fulfill multiple requests introduce high degrees of uncertainty, creating challenges for ride-pooling operators. This study aims to expand the current state-of-the-art Approximate Dynamic Programming (ADP) approach for ride-pooling services, introduce key extensions, and perform a comparative analysis with the Neural Approximate Dynamic Programming (NeurADP) approach to optimize the efficiency and effectiveness of these services. Specifically, we develop an ADP approach that incorporates three important problem specifications: (i) pick-up and drop-off deadlines, (ii) vehicle rebalancing, and (iii) allowing more than two passengers in a vehicle. We conduct a detailed numerical study with the New York City taxi-cab dataset and a novel dataset of taxi-cab requests collected in the city of Chicago. We also provide a sensitivity analysis on key model parameters such as wait and delay times, passenger group sizes, and vehicle capacity, along with the investigation of the effects of vehicle rebalancing. Our comparative analysis highlights the strengths and limitations of both ADP and NeurADP methodologies. Network density and road directionality are found to significantly impact the performance. NeurADP is found to be more efficient in learning value function approximations for larger and more complex problem settings than the ADP approach. However, for smaller settings, ADP is shown to outperform NeurADP.
\end{abstract}
\end{frontmatter}
\vspace{-8pt}
\section{Introduction}
Prominent corporations such as Uber and Lyft provide on-demand ride-pooling services which facilitate the combination of passenger requests and enable groups of travelers to share transportation to their respective destinations. In addition to mitigating emissions and alleviating traffic congestion, ride-pooling services such as UberPool have also been shown to confer significant benefits to both passengers and corporations. According to a study conducted in 2018, UberPool has reduced over 36,537 tonnes of CO2 emissions since its launch in India in 2014, underscoring the potential of such services to contribute to environmental sustainability efforts \citep{agarwal2018}. Passengers who use ride-pooling services also stand to benefit from reduced travel costs as a result of the shared mobility model, while corporations may enjoy higher revenues due to the ability to fulfill more requests with fewer vehicles. These benefits have contributed to the increasing popularity of ride-pooling services, and despite a brief halt in services during the COVID-19 pandemic, market projections indicate that the ride-sharing industry is expected to reach a valuation of \$242.73 billion USD by 2028 \citep{heath2016, fortune2021}.

As ride-pooling services continue to grow in popularity, it is increasingly important to prioritize efficiency in operations. The main objective of ride-pooling operations is to maximize the number of passengers served while minimizing wait and delay times. However, the nature of ride-pooling services introduces a high degree of uncertainty, stemming from factors such as the timing and volume of passenger requests, pick-up and drop-off locations, available capacity in vehicles, and the trajectory needed to fulfill multiple requests. These uncertainties pose significant challenges for ride-pooling operators, who must navigate them in order to provide efficient and effective services to passengers. To address these challenges, this paper makes several contributions. 

First, we expand upon the current state-of-the-art Approximate Dynamic Programming (ADP) approach taken by \citet{yu2019integrated} to solve the ride pooling problem as a multistage stochastic program. More specifically, we introduce three novel extensions: (i) pick-up and drop-off deadlines for passenger requests to enable accurate modeling of real-world time constraints faced by ride-pooling services, (ii) vehicle rebalancing to improve operational efficiency through the relocation of unassigned vehicles, and (iii) allowance of more than two passenger groups in the same vehicle, which helps to obtain a more realistic model in practice. Additionally, we perform sensitivity analysis on several key model parameters, including wait and delay times, passenger group sizes, and vehicle capacity, as well as investigate the effects of enabling vehicle rebalancing. By exploring the impact of these parameters on the efficiency and effectiveness of ride-pooling services, we gain insights into how best to optimize the provision of these services for the benefit of both passengers and ride-pooling companies. Furthermore, this paper introduces algorithmic enhancements to the ADP approach for the Ride-Pool Matching Problem (RMP) in the form of hierarchical aggregation (HA) and bias-adjusted Kalman filter (BAKF) step size, which lead to improved computational efficiency and faster convergence speed.

Second, we provide a comparative analysis between our traditional ADP approach and the Neural Approximate Dynamic Programming (NeurADP) approach introduced by \citet{shah2020neural}, which employs neural networks to perform value function approximation. While previous studies have compared NeurADP-based approaches to myopic solutions, our work provides the first comparative analysis with ADP and NeurADP solution methodologies. Such thorough comparative analysis is critical for future practitioners to have a comprehensive understanding of the strengths, weaknesses, trade-offs, and limitations of the two solution methodologies. To facilitate this comparison, we design a hybrid problem setting which entails one-to-one vehicle-to-request matching at each time step, enabling the development of a more flexible action selection. Additionally, this setting ensures immediate pick-up of newly matched passenger requests, reducing the computational complexity of generating feasible actions.

Lastly, in addition to the commonly used New York City taxi-cab dataset, we introduce and conduct experiments on a novel dataset of taxi-cab requests collected in the city of Chicago. This dataset provides valuable information on the daily distribution of ride requests and the frequency of pick-up and drop-off locations. We then map the distribution of requests onto a grid system and analyze how the directionality and denseness of the network impact the performance of the ADP and baseline solution methodologies.  

The rest of the paper is structured as follows. Related literature is presented in Section~\ref{literaturereview}; a formal problem description and problem setting are given in Section~\ref{problemdescription} and Section~\ref{problemformulation}, respectively. The ADP solution methodology is described in Section~\ref{solutionmethodology}. Section~\ref{experimentalsetup} provides details about the datasets and the benchmark policies. Finally, Section~\ref{results} presents the results of the computational experiments, and Section~\ref{conclusion} concludes the paper by summarizing the research findings and suggesting directions for future research.

\section{Literature Review} \label{literaturereview}
The optimization problems of Taxi-on-Demand (ToD) and RMP are significant focus areas in transportation planning and ride-sharing services. ToD aims to dispatch taxis to customer requests efficiently, while RMP focuses on optimizing the matching of multiple riders with drivers to minimize the required number of vehicles. Unlike ToD, which can handle only one passenger group at a time, RMP can pool multiple passenger groups together at the same time, leading to higher utilization of available vehicles. In both cases, efficient algorithms are required to balance various factors such as customer wait times, vehicle idle times, operational costs, and environmental impacts. To address these problems, a variety of approaches have been proposed, including heuristic algorithms, machine learning models, combinatorial optimization models, reinforcement learning algorithms, and market-based mechanisms. In Table~\ref{table:references}, the relevant studies for both problems are presented. This table is comprised of six indicators that provide information about the problem setting and solution methodology. These are: ``Solution Technique'', which describes the approach used to solve the problem, ``Pooling'', which indicates whether the problem is designed for the RMP, ``Non-Myopic'', which indicates whether a myopic solution technique is employed, ``Batching'', which indicates whether passenger requests are considered in batches or sequentially, ``High Group Capacity'', which indicates whether the problem permits more than two distinct groups of passenger requests inside a vehicle at any given time, and ``Rebalancing'', which indicates whether idle vehicles are relocated.

\setlength{\tabcolsep}{2.5pt} 
\renewcommand{\arraystretch}{1.3} 
\begin{table}[!ht]
\centering
\caption{Summary of relevant studies. (CH: Combinatorial Heuristic, RL: Reinforcement Learning, DP: Dynamic Programming, ADP: Approximate Dynamic Programming, NeurADP: Neural Approximate Dynamic Programming)}\label{table:references}
\resizebox{0.99\textwidth}{!}{
\begin{tabular}{P{0.29\textwidth}C{0.15\textwidth}C{0.15\textwidth}C{0.15\textwidth}C{0.15\textwidth}C{0.15\textwidth}C{0.15\textwidth}} 
\toprule
\textbf{Study} & \textbf{Solution Technique} & \textbf{Pooling} & \textbf{Non-Myopic} & \textbf{Batching} & \textbf{High Group Capacity} & \textbf{Rebalancing} \\ 
\midrule
\citet{alonso2017demand}      & CH & \checkmark                           &                                         & \checkmark                            & \checkmark                                 & \checkmark                               \\
\citet{zhang2017taxi}         & CH &                                      &                                         &                                       &                                            &                                          \\
\citet{chen2018ptrider}       & CH & \checkmark                           &                                         &                                       & \checkmark                                 &                                          \\
\citet{lin2018efficient}      & RL &                                      & \checkmark                              &                                       &                                            & \checkmark                               \\
\citet{sayarshad2018scalable} & DP &    \checkmark                                  & \checkmark                              &                                       &                                            & \checkmark                               \\
\citet{tong2018unified}       & DP & \checkmark                           & \checkmark                              &                                       & \checkmark                                 &                                          \\
\citet{zheng2018order}        & CH & \checkmark                           &                                         &                                       & \checkmark                                 &                                          \\
\citet{lowalekar2019zac}      & CH & \checkmark                           &                                         & \checkmark                            & \checkmark                                 &                                          \\
\citet{yu2019integrated}      & ADP                                             & \checkmark                           & \checkmark                              & \checkmark                            &                                            &                                          \\
\citet{al2020approximate}     & ADP                                             &                                      & \checkmark                              & \checkmark                            &                                            & \checkmark                               \\
\citet{ke2020learning}        & RL &                                      & \checkmark                              &                                       &                                            &                                          \\
\citet{shah2020neural}        & NeurADP                                         & \checkmark                           & \checkmark                              & \checkmark                            & \checkmark                                 & \checkmark                               \\
\citet{bose2021conditional}   & NeurADP                                         & \checkmark                           & \checkmark                              & \checkmark                            & \checkmark                                 & \checkmark                               \\
\citet{kumar2021improving}    & NeurADP                                         & \checkmark                           & \checkmark                              & \checkmark                            & \checkmark                                 & \checkmark                               \\
\citet{qin2021optimizing}     & RL &                                      & \checkmark                              & \checkmark                            &                                            &                                          \\
\citet{hao2022hierarchical}   & NeurADP                                         & \checkmark                           & \checkmark                              & \checkmark                            & \checkmark                                 & \checkmark                               \\ 
\midrule
\textbf{Our Work}               & ADP                                             & \checkmark                           & \checkmark                              & \checkmark                            & \checkmark                                 & \checkmark                               \\ 
\bottomrule
\end{tabular}
}
\end{table}

Numerous studies focused on tackling the problem setting induced in ToD services. While some have proposed myopic solution strategies, such as the one presented by \citet{zhang2017taxi}, which introduced a taxi dispatch system that handles multiple bookings and predicts user destinations using a Bayesian framework based on travel histories, the majority of studies concentrate on non-myopic solutions. For instance, both \citet{ke2020learning} and \citet{qin2021optimizing} proposed reinforcement learning-based solutions that concentrate on the effect of delayed matching between vehicles and passenger requests. In particular, \citet{ke2020learning} proposed a two-stage framework for matching passengers and drivers in ride-sourcing systems that considers various metrics and enhances matching effectiveness. On the other hand, \citet{qin2021optimizing} proposed reinforcement learning algorithms for delayed matching in ride-hailing systems that balances wait time penalties and efficiency, with the added capability of being incorporated as a lookup table for adaptive decision making. Subsequently, \citet{lin2018efficient} focused on the assignment of agents to different zones and proposed a reinforcement learning-based approach for fleet management in large-scale ride-sharing platforms using a multi-agent framework and two coordination algorithms. Finally, \citet{al2020approximate} employed an ADP technique for a ride-sharing system that runs on a self-driving electric vehicle fleet. Their method involved imposing monotonicity on the value functions and incorporating hierarchical aggregation.

On the other hand, the RMP presents a more intricate problem setting in which pre-occupied vehicles are regularly rerouted to pick up and drop off newly assigned passengers. As a result, most of the earlier studies on RMP considered myopic approaches. \citet{alonso2017demand} introduced a reactive anytime myopic method for assigning passenger requests to a fleet of vehicles with varying capacities and investigated the trade-off between fleet size, capacity, waiting time, travel delay, and operational costs. \citet{chen2018ptrider} presented PTRider, a price-and-time-aware ride-sharing system that offers travellers more options by considering pick-up time and price. It builds indexes on the road network and vehicles to match requests in real time and uses a real-life dataset to demonstrate its effectiveness. \citet{zheng2018order} considered the price of the orders as an important reference to the platform's profit. They formulated constrained optimization models, developed approximation methods, and evaluated the effectiveness and efficiency of these methods through experiments on a simulation framework using real ride-sharing data. Finally, \citet{lowalekar2019zac} proposed ``ZAC'', a ride-sharing method which clusters locations into zones and utilizes paths of zones instead of individual location points. The approach generates multiple zone paths to represent various trips and assigns available vehicles to optimize the fulfillment of ride requests.

To address the issue of myopic decision-making in the RMP, various studies have proposed dynamic programming (DP)-based approaches. These methods are able to handle the problem's combinatorial complexity and incorporate future information into the decision-making process. For instance, \citet{sayarshad2018scalable} suggested a competitive on-demand mobility system which utilizes a multi-server queue and a pricing policy based on operational and willingness-to-pay factors. In contrast, \citet{tong2018unified} proposed a DP-based algorithm for shared mobility route planning. This method enhances operational efficiency by introducing a multi-objective function. However, DP-based approaches have limitations when it comes to solving large-scale RMPs, and as such, these approaches may not be adequate for larger-scale versions of the problem. As a remedy, Approximate Dynamic Programming (ADP) has been employed to reduce the overall problem complexity by performing value function approximation over post-decision states. On this front, \citet{yu2019integrated} proposed an ADP-based solution for solving the ride-pooling problem, where two passenger groups share rides simultaneously. In order to simplify the complexity of the problem, the solution utilizes a decomposition heuristic which breaks down the decision horizon into distinct stages and spatially partitions the map into sub-regions.

Another novel approach applied to the RMP is Neural Approximate Dynamic Programming (NeurADP), which was introduced by \citet{shah2020neural}. NeurADP utilizes a neural network-based approximate value function to handle the combinatorial complexity of passenger requests. \citet{hao2022hierarchical} built upon this approach by proposing HIerarchical ValuE decompoSition (HIVES), which considered the future impact of current matches and the impact of other agents on the future value of a vehicle via the use of mixing networks and spatial clustering. \citet{bose2021conditional} proposed a novel mechanism called Conditional Expectation-based Value Decomposition (CEVD) that captures dependencies on other agents' actions without increasing the complexity of training or decision making. Finally, \citet{kumar2021improving} proposed an online approach which incorporated fairness components into the value function during matching optimization, which eliminates the need for retraining and can be adapted to use any existing value function approximator. \citet{wang2023optimization} studied an extension of the traditional ride-sharing problem where passenger transfer between vehicles is allowed at transfer stations. They employed a joint decision framework that combines deep reinforcement learning and integer programming to solve their problem. 
\citet{hua2022optimality} proposed mixed-integer models for the dynamic shared-taxi problem where they consider rearrangement of the requests to assess the trade-off between system-wide profit and service quality.
\citet{najmi2017novel} took a different approach to the real-time ride-sharing problem which iteratively solves a matching problem in a rolling horizon approach.

Our work builds upon the current state-of-the-art ADP approach for the RMP introduced by \citet{yu2019integrated}. Specifically, we present a hybrid one-to-one vehicle-to-request setting of the RMP and propose three key extensions. First, we extend the limited two-passenger group setting considered by \citet{yu2019integrated} and allow for the maximum passenger group size to be a parameter that can be increased beyond two groups. Second, we introduce pick-up and drop-off deadlines for passenger requests to incorporate real-world time constraints faced by ride-pooling services. Third, we introduce a vehicle rebalancing strategy that allows idle vehicles to relocate to areas of higher demand, thereby improving operational efficiency and enabling more accurate modeling of real-world scenarios. We also introduce algorithmic enhancements to the ADP framework in the form of hierarchical aggregation (HA) and bias-adjusted Kalman filter (BAKF) step size. Next, we connect the ADP and NeurADP works in the literature by providing a comparative analysis of the two methods and conducting numerical experiments to the backdrop of the standard New York City yellow taxicab dataset, as well as a novel Chicago taxicab dataset that we introduce. Finally, we map the distribution of requests onto a grid and analyze how the directionality and denseness of the network impact the performance of different solution methods for RMP.

\section{Problem Description} \label{problemdescription}
We present a dynamic ride-pooling model which seeks to match vehicles with incoming batches of passenger requests. The model considers the spatial and temporal demand patterns of the requests, which arrive dynamically over a multi-period decision horizon. Passenger requests, which may correspond to multiple passengers with identical pickup and drop-off locations, are generated via a stochastic process and must be picked up and dropped off within corresponding deadlines. 
Furthermore, there exists a predetermined set of vehicles available during the planning horizon, subject to capacity limitations. 

We make three main assumptions in our model. First, we assume that the matching of vehicles and requests at every time step is done in a one-to-one manner. That is, an available vehicle can be matched to at most one request at each decision epoch. This is a reasonable assumption given a sufficiently short decision epoch, which we set to 60 seconds in our experiments. Our second assumption is predicated upon the rerouting required to determine the most time-efficient sequence for picking up and dropping off passenger requests. We assume that vehicles must promptly pick up new requests upon being matched, before dropping off any previously assigned passenger groups or being matched with any new ones. This assumption, as was also made by \citet{yu2019integrated}, is a simplification which reduces computational effort and combinatorial complexity during the matching process. Furthermore, it is reasonable to assume that this approach will not result in exceeding the drop-off deadlines and tight pick-up deadlines will not impose significant detours for onboard passengers. Finally, in line with previous RMP research, we assume that any requests which are not matched within the period of their arrival exit the system. This assumption is reasonable since customers are generally unwilling to wait for an extended amount of time to receive a confirmation on whether their requests have been accepted by the system or not. 

The main objective of our model is to maximize the total number of passenger requests served within the decision horizon. To achieve this goal, our assignment decisions take into account future uncertainties and the potential downstream costs of current decisions. To address the decision complexity inherent in this approach, we formulate a Markov Decision Process (MDP) model and develop an ADP framework that enables effective real-time decision making under uncertainty.

\section{Problem Formulation} \label{problemformulation}
We formulate an MDP model for the dynamic ride-pooling problem considering a problem setting wherein vehicles travel throughout the city to fulfill incoming passenger requests. To achieve this, we divide the finite planning horizon into discrete time intervals, each of length \TimeIntervals~(e.g., 60 seconds), and assume decisions are made at the beginning of each interval, while exogenous information is observed throughout. At the end of each interval, the state of the system is updated based on the decisions made and exogenous information observed. The set of decision-making epochs is represented by $\Time:=\{1,\ldots,T\}$. We define the network where vehicles and passengers interact as $\GraphNetwork=(\Locations,\Edges)$, where $\Locations$ represents the set of street intersections, and $\Edges$ represents the adjacency of these intersections. The weights assigned to the connections between intersections correspond to the travel time of the road segment. Given two intersection locations $\ell$, $\ell'$ $\in \Locations$, we denote the travel time by $\TravelTime(\ell,\ell')$. It is assumed that vehicles can only collect and drop off passengers at these intersections. We establish uniform definitions for \MaxWait~ and \MaxDelay~ across all requests, representing the maximum acceptable wait and delay times for a given request, respectively. More specifically, \MaxWait~ pertains to the maximum duration within which a vehicle must pick-up a passenger once the request has entered the system and been matched with an agent. On the other hand, \MaxDelay~ signifies the maximum additional duration beyond the original travel time within which the passenger must be dropped off after being picked up. Any unused wait time is added to the delay time. As an illustration, consider a scenario where the wait and delay times are set to $60$ seconds and the original travel time for a request is $40$ seconds. In this case, the vehicle must pick-up the passenger(s) within $60$ seconds of the request being made. If the pick-up is completed within $15$ seconds, the total duration available from that point onward to complete the drop-off would be $45 + 40 + 60 = 145$ seconds. 

We define \MaxVehicleGroups~ to be the maximum number of unique passenger groups, each having their own designated pick-up and drop-off points, that can be accommodated by a single vehicle at any given time. It is noteworthy that all passengers within a group share the same pickup and drop-off points. On the other hand, we define \MaxVehicleCapacity~ to relate to the restriction placed on the number of total individuals, irrespective of the number of groups, a vehicle can fit at any time. For instance, suppose a vehicle has a group size limit of $\MaxVehicleGroups=2$ and a capacity size limit of $\MaxVehicleCapacity=5$. In this case, the vehicle can accommodate up to two distinct requests at any given time, with only 5 seats available for passengers to occupy. As an illustration, the first request may consist of 2 passengers, while the second request may have 3 passengers. Both \MaxVehicleGroups~ and \MaxVehicleCapacity~ are set uniformly for all vehicles in the system. Finally, we establish the set $\RebalancingPointsSet$ of rebalancing points within the city. Rebalancing refers to the ability of an unassigned vehicle to reposition itself to a new position. In practice, we partition the map into distinct zones and identify high-demand locations within each zone based upon historical data. At every time-step, a vehicle has the flexibility to move towards any of these high-demand locations. By doing so, the vehicle can place itself in areas with a greater potential demand for service, which in turn can reduce the time required to pick up passengers and assist in the main objective of maximizing the number of requests served. In the following subsections, we provide further details on the various components that make up the MDP model, including the state variables, decision variables, cost function, exogenous information, and transition function.
\subsection{State Variables}
The state of the system is characterized by online passenger requests and vehicles. The state of an individual vehicle can be represented as a two-dimensional attribute vector denoted by $\IndividualVehicle=(\IndividualVehicle_\VehicleLocation,\IndividualVehicle_\VehicleRequests)\in \VehiclesSet$. Here, $\IndividualVehicle_\VehicleLocation \in \Locations$ represents the current position of the vehicle, while $\IndividualVehicle_\VehicleRequests$ is an ordered list of passenger requests assigned to the vehicle. The requests are arranged so as to maximize delivery efficiency, and the list is re-ordered whenever a new request is assigned. Subsequently, $\IndividualVehicle_\VehicleRequests$ contains important details about the vehicle's current tasks, such as the requests it is currently assigned to, pending pickups, deadlines, and the number of passengers for each assigned request. Subsequently, the state of a passenger request is expressed by a five-dimensional attribute vector denoted by $\IndividualRequest=(\IndividualRequest_\RequestOrigins, \IndividualRequest_\RequestDestinations, \IndividualRequest_\RequestCapacity, \IndividualRequest_\RequestDeadline, \IndividualRequest_\IndividualVehiclePickedUp) \in \RequestsSet$. Here, $(\IndividualRequest_\RequestOrigins, \IndividualRequest_\RequestDestinations)\in(\Locations \times \Locations)$ refer to the pick-up and drop-off locations of the request, respectively, while $\IndividualRequest_\RequestCapacity$ denotes the number of passengers who are associated with the request. The $\IndividualRequest_\RequestDeadline$ component represents the deadline by which the passenger request must be completed. When a passenger request arrives between decision epochs $t-1$ and $t$, its deadline attribute $\IndividualRequest_\RequestDeadline$ is determined at the beginning of epoch $t$ using the equation
\begin{equation}
    \IndividualRequest_\RequestDeadline = \CurrentTime + \MaxWait + \TravelTime(\IndividualRequest_\RequestOrigins, \IndividualRequest_\RequestDestinations) + \MaxDelay \label{eq:PassDeadline}
\end{equation}
Finally, $\IndividualRequest_\IndividualVehiclePickedUp$ is a binary value denoting whether or not the passenger request has been picked up by a vehicle. This value is set to 0 by default until the request is matched with a vehicle and the corresponding passenger is picked up. Then, let $\StateNote_{\CurrentTime \IndividualVehicle}^{\texttt{Vehicle}}$ and $\StateNote_{\CurrentTime \IndividualRequest}^{\texttt{Request}}$ be the number of vehicles with attribute vector $\IndividualVehicle \in \VehiclesSet$ and requests with attribute $\IndividualRequest \in \RequestsSet$ at decision epoch $\CurrentTime \in \Time$, respectively. The system state vector at $\CurrentTime \in \Time$, before decisions are made, is $\StateNote_{\CurrentTime} := (\StateNote_{\CurrentTime}^{\texttt{Vehicle}},\StateNote_{\CurrentTime}^{\texttt{Request}})$, which consists of the vehicles state vector, $\StateNote_{\CurrentTime}^{\texttt{Vehicle}} = (\StateNote_{\CurrentTime \IndividualVehicle}^{\texttt{Vehicle}})_{\IndividualVehicle \in \VehiclesSet}$ and the passenger requests state vector, $\StateNote_{\CurrentTime}^{\texttt{Request}} = (\StateNote_{\CurrentTime \IndividualRequest}^{\texttt{Request}})_{\IndividualRequest \in \RequestsSet}$.
\subsection{Decision Variables}
For each decision epoch $\CurrentTime \in \Time$, we determine a matching between available vehicles and a batch of incoming passenger requests based on the current state of the system. To determine whether it is possible to match a vehicle with state $\IndividualVehicle$ to a new request with state $\IndividualRequest$, we consider the feasibility of adding $\IndividualRequest$ to $\IndividualVehicle$'s ordered list of passenger requests, represented by $\IndividualVehicle'_\VehicleRequests$, such that $\IndividualVehicle'_\VehicleRequests = \IndividualVehicle_\VehicleRequests \cup \{\IndividualRequest\}$. This decision is subject to the following constraints:
\begin{subequations}
\label{m:FeasibilityConsts}
    \begin{alignat}{2}
    & \sum_{\IndividualRequest' \in \IndividualVehicle'_\VehicleRequests}
    \IndividualRequest'_\RequestCapacity \leq \MaxVehicleCapacity  \label{eq: CapacityConst} \\
    & \big| \IndividualVehicle'_\VehicleRequests  \big| \leq \MaxVehicleGroups
    \label{eq: GroupsConst} \\
    & \sum_{\IndividualRequest' \in \IndividualVehicle'_\VehicleRequests}
    \IndividualRequest'_\IndividualVehiclePickedUp = 1   
    \label{eq: PickupConst} \\
    & \TravelTime(\IndividualVehicle_\VehicleLocation, \IndividualRequest_\RequestOrigins) \leq \MaxWait
    \label{eq: WaitConst} \\
    & \IndividualRequest'_\RequestDeadline \leq \IndividualRequest'_\IndividualVehicleDropOff && \quad \quad  \quad \forall \IndividualRequest' \in \IndividualVehicle'_\VehicleRequests
    \label{eq: DelayConst}
    \end{alignat}
\end{subequations}
The constraint expressed in \eqref{eq: CapacityConst} guarantees that the vehicle, which may already have passengers from prior requests, has sufficient capacity to take on additional passengers for the new request. The subsequent constraints, \eqref{eq: GroupsConst} and \eqref{eq: PickupConst}, ensure that the addition of a new passenger group does not exceed the allowable number of unique passenger groups, and that the vehicle is not in the process of picking up any previously assigned requests. Finally, constraints \eqref{eq: WaitConst} and \eqref{eq: DelayConst} ensure that the vehicle can pick up the new request before its wait deadline and drop off all assigned requests prior to their respective deadlines. Here, $\IndividualRequest'_\IndividualVehicleDropOff$ represents the time of drop-off for each assigned request, which takes into account both the current position of the vehicle and the revised route it will take based on the updated ordering of passenger requests in $\IndividualVehicle'_\VehicleRequests$. We can define the feasible set of matches between vehicles and requests at time $\CurrentTime$ as:
\begin{equation}
\FeasibilitySet_\CurrentTime = \{(\IndividualVehicle,\IndividualRequest) \in (\VehiclesSet_\CurrentTime \times \RequestsSet_\CurrentTime)\ : (\eqref{eq: CapacityConst}-\eqref{eq: DelayConst})\}
\label{eq:FeasibilitySet}
\end{equation}
Additionally, for each un-assigned vehicle, we can define a set of feasible relocation matchings as follows:
\begin{equation}
\RebalancingFeasibilitySet_\CurrentTime = \{(\IndividualVehicle,\RebalancingPoint) \in (\VehiclesSet_\CurrentTime,\RebalancingPointsSet) : \IndividualVehicle_\VehicleRequests = \emptyset \}
\label{eq:RelocationSet}
\end{equation}
Finally, at every time step, each vehicle has a null action defined by $\NullAction$ which signifies the agent continuing on their given path. If the vehicle is unassigned, it will remain stationary at its current location. On the other hand, if it is assigned to pick up and deliver passenger requests, it will proceed along its path towards the designated locations. Therefore, the decision set at time $\CurrentTime$ is represented by $\MovementSet_\CurrentTime = \RequestsSet_\CurrentTime ~\cup \RebalancingPointsSet ~\cup {\NullAction}$. Subsequently, at each decision epoch $\CurrentTime$, the vehicles and passenger requests can be used to formulate a linear program which matches vehicles with passengers and determines potential rebalancing strategies:
\begin{subequations}
\label{m:LinearSystem}
    \begin{alignat}{2}
    & \sum_{\IndividualRequest \in \RequestsSet_\CurrentTime: (\IndividualVehicle, \IndividualRequest) \in \FeasibilitySet_\CurrentTime} \MatchingVariable_{\CurrentTime\IndividualVehicle\IndividualRequest} +  \sum_{(\IndividualVehicle,\RebalancingPoint) \in \RebalancingFeasibilitySet_\CurrentTime} \MatchingVariable_{\CurrentTime\IndividualVehicle\RebalancingPoint} + \MatchingVariable_{\CurrentTime \IndividualVehicle \NullAction} = \StateNote_{\CurrentTime \IndividualVehicle}^{\texttt{Vehicle}} && \quad \quad \quad \forall \IndividualVehicle \in \VehiclesSet_\CurrentTime \label{eq: VehicleFlow} \\
    & \sum_{\IndividualVehicle \in \VehiclesSet_\CurrentTime: (\IndividualVehicle, \IndividualRequest) \in \FeasibilitySet_\CurrentTime} \MatchingVariable_{\CurrentTime\IndividualVehicle\IndividualRequest} \leq \StateNote_{\CurrentTime \IndividualRequest}^{\texttt{Request}} && \quad \quad \quad \forall \IndividualRequest \in \RequestsSet_\CurrentTime \label{eq: PassengerFlow} \\
    & \MatchingVariable_{\CurrentTime\IndividualVehicle\MovementPoint} \geq 0 && \quad \quad \quad \forall \IndividualVehicle \in \VehiclesSet_\CurrentTime, \MovementPoint \in \MovementSet_\CurrentTime \label{eq: VariableDefinition}
    \end{alignat}
\end{subequations}

Constraints~\eqref{eq: VehicleFlow} and \eqref{eq: PassengerFlow} represent the flow conservation for vehicles and passengers, respectively. More specifically, \eqref{eq: VehicleFlow} ensures that the number of actions assigned to a vehicle with a certain attribute is equal to the number of vehicles with said attribute type. Similarly, \eqref{eq: PassengerFlow} asserts that the number of drivers assigned to passengers with a specific attribute does not exceed the total number of passengers with that attribute. Finally, \eqref{eq: VariableDefinition} ensures that all decision variables are non-negative. Given that the coefficient matrix formed by constraints \eqref{eq: VehicleFlow} and \eqref{eq: PassengerFlow} is totally unimodular and the right-hand side vector is integral, integral solutions can be obtained at all extreme points of the linear program~\citep{wolsey1999integer}. Allowing $\DecisionTuple^\CurrentTime \in \DecisionSet_\CurrentTime$ to represent the set of actions taken at time $\CurrentTime$, we can define the total reward collected at each time step $\CurrentTime \in \Time$ by
\begin{equation}
    \Reward (\DecisionTuple^\CurrentTime)=\sum_{(\IndividualVehicle, \IndividualRequest) \in \FeasibilitySet_\CurrentTime} \MatchingVariable_{\CurrentTime\IndividualVehicle\IndividualRequest}~,
    \label{eq: RewardFunction}
\end{equation}
which characterizes the number of matched passenger requests.


\subsection{Exogenous Information}
Exogenous information takes the form of incoming passenger requests, which arrive at every time step during the decision horizon. It is defined by $\ExogenousInformation_\CurrentTime = (\ExogenousInformation_{\CurrentTime\IndividualRequest}^\texttt{Request})_{\IndividualRequest \in \RequestsSet_\CurrentTime}$, in which $\ExogenousInformation_{\CurrentTime\IndividualRequest}^\texttt{Request}$ represents the number of passenger requests with attribute type $\IndividualRequest$ arriving between decision epochs $\CurrentTime$ and $\CurrentTime + 1$.

\subsection{Transition Function}

The transition function, which defines the evolution of the state of the system from time $\CurrentTime$ to $\CurrentTime + 1$, depends upon the arrival of passenger requests and the decision tuple $\DecisionTuple^\CurrentTime \in \DecisionSet^\CurrentTime$. We are able to break down this evolution into two parts through the introduction of a post-decision state \citep{powell2007approximate}, which defines the state of the system immediately following a decision but prior to the arrival of exogenous information in the next time step. The first transition \eqref{eq: FirstTrans} is to the post-decision state, $\PostDecisionState_\CurrentTime$, via action $\DecisionTuple^\CurrentTime$, while the second transition \eqref{eq: SecondTrans} is from the post-decision state to the next state based upon the arrival of exogenous information $\ExogenousInformation_\CurrentTime$.
\begin{subequations}
\label{m:Transitions}
    \begin{alignat}{2}
    & \PostDecisionState_\CurrentTime=\FirstTransition(\StateNote_\CurrentTime,\DecisionTuple^\CurrentTime)  \label{eq: FirstTrans} \\
    & \StateNote_{\CurrentTime+1}=\SecondTransition(\PostDecisionState_\CurrentTime,\ExogenousInformation_\CurrentTime)  \label{eq: SecondTrans}
    \end{alignat}
\end{subequations}
As previously mentioned, it is assumed that unassigned passenger requests leave the system at the end of each decision epoch, and thus, the state of passenger requests at time $\CurrentTime + 1$ is defined by:
\begin{equation}
    \label{eq: RequestTransition}
    \StateNote_{{\CurrentTime + 1},\IndividualRequest}^{\texttt{Request}} = \ExogenousInformation_{\CurrentTime \IndividualRequest} \quad \quad \quad \forall \IndividualRequest \in \RequestsSet_\CurrentTime
\end{equation}
Additionally, the state of vehicles at time $\CurrentTime + 1$ is defined by:
\begin{equation}
    \label{eq: VehicleTransition}
    \StateNote_{{\CurrentTime + 1},\IndividualVehicle}^{\texttt{Vehicle}} = \StateNote^{\VehiclePost}_{\CurrentTime\IndividualVehicle} \quad \quad \quad \forall \IndividualVehicle \in \VehiclesSet_\CurrentTime
\end{equation}
such that $\StateNote^{\VehiclePost}_{\CurrentTime\IndividualVehicle}$ denotes the state of vehicle $\IndividualVehicle$ after being assigned an action and simulating its motion forward in time by $\TimeIntervals$, prior to the arrival of new exogenous information. For an idle vehicle, if the assigned action is to remain put, its state remains the same. On the other hand, if the assigned action is to rebalance to a new location, the vehicle's current position is updated according to its motion towards the new location. If a passenger request is assigned to an idle or non-idle vehicle, the vehicle's ordered list of requests is updated, and the vehicle moves to pick up and drop off the updated set of requests.
\subsection{Optimal Policy}
The objective in our ride-pooling problem is to maximize the expected number of requests served over the operation horizon:
\begin{equation}\label{eq:ADPObj}
    \max_{\pi \in \Pi}\mathbb{E}_{\ExogenousInformation=(\ExogenousInformation_1, \ldots, \ExogenousInformation_{T})} \left[ \sum_{\CurrentTime \in \Time} \Reward \left( \DecisionSet_\CurrentTime^\pi(\StateNote_\CurrentTime^\pi(\ExogenousInformation)) \right) \big|\StateNote_1 \right]
\end{equation} 
By solving the Equation~\eqref{eq:ADPObj}, we can determine a policy $\pi$ from a set of feasible policies $\Pi$ which yields maximum reward when its suggested actions $\DecisionSet_\CurrentTime^\pi(\StateNote_\CurrentTime)$ are implemented sequentially at realized states. The realized states are such that:
\begin{subequations}
\label{m:Transitions2}
    \begin{alignat}{2}
    & \StateNote_1^\pi(\ExogenousInformation)=\StateNote_{1}  \label{eq: OptOne} \\
    & \StateNote_{\CurrentTime+1}^\pi(\ExogenousInformation)=\SecondTransition(\FirstTransition(\StateNote_\CurrentTime^\pi, \DecisionSet_\CurrentTime^\pi(\StateNote_\CurrentTime^\pi)),\ExogenousInformation_{\CurrentTime})  \label{eq: OptTwo}
    \end{alignat}
\end{subequations}
for $\CurrentTime \in \{1, \ldots, T\}$, where $\StateNote_1$ represents the initial state in the decision horizon, which contains initialized vehicles as well as prior requests. The future reward is obtained based upon the expectation with respect to the underlying stochastic process described by $\ExogenousInformation$. The actions taken and states visited during a decision epoch are dependent only on the revealed random variables up to that point, while the total reward depends on the realization of the complete vector $\ExogenousInformation$. An optimal policy can be obtained through backward DP \citep{powell2007approximate} via solving the Bellman optimality equation:
\begin{equation}
V_\CurrentTime(\StateNote_\CurrentTime)=\max_{\DecisionTuple^\CurrentTime \in \DecisionSet^\CurrentTime(\StateNote_\CurrentTime)}\Reward(\DecisionTuple^\CurrentTime) + \mathbb{E}_{\ExogenousInformation_{\CurrentTime}}[V_{\CurrentTime + 1}(\StateNote_{\CurrentTime + 1})]
    \label{eq: Bellman}
\end{equation}
where $\StateNote_{\CurrentTime + 1} = \SecondTransition(\FirstTransition(\StateNote_\CurrentTime,\DecisionTuple^\CurrentTime),\ExogenousInformation_\CurrentTime)$.   
That is, in order to update value function of a particular state, denoted by $V_t(S_t)$, a procedure is utilized which involves working backwards in time from the final epoch, $T$, while taking into account the rewards associated with making the optimal set of decisions, as well as the probability of transitioning from one state to another. This recursive process continues until the first stage, $t=1$, is reached. However, this approach becomes impractical for even small instances due to the need for enumerating all possible outcomes and actions. Therefore, we propose an ADP policy that employs value function approximation for post-decision states as a more practical alternative to calculating the exact value function for each state.
\section{Solution Methodology} \label{solutionmethodology}
The MDP formulation for the ride-pooling problem suffers from ``three curses of dimensionality'', resulting from intractability in the state, action, and outcome space~\citep{powell2007approximate}. To address this challenge, we adopt an ADP approach which combines reinforcement learning and mathematical programming techniques, and entails an offline learning process of value functions through simulation of the problem's Markov decision process. In this section, we first provide a description of our post-decision states and value function approximations, then we introduce our ADP policy and forward pass algorithm. Finally, we provide a set of algorithmic enhancements to our ADP framework.

\subsection{Post-decision State}
In order to solve the Bellman optimality equation provided in~\eqref{eq: Bellman} for state $\StateNote_\CurrentTime$ via dynamic programming, the anticipated downstream reward needs to be computed. This necessitates multiplying the value function of every conceivable outcome of $\StateNote_{\CurrentTime + 1}$ by the probability associated with arriving at $\StateNote_{\CurrentTime + 1}$, which is determined by the exogenous information seen $\ExogenousInformation_\CurrentTime$. Nevertheless, for large-scale problems such as the RMP, the outcome space may become excessively large, leading to first curse of dimensionality. To address this problem, \citet{powell2007approximate} put forward the notion of the post-decision state. This technique allows us to divide the dynamic programming recursion equation into two parts:
\begin{subequations}
\label{m:BellmanBreakup}
    \begin{alignat}{2}
    & V_\CurrentTime(\StateNote_\CurrentTime)=\min_{\DecisionTuple^\CurrentTime \in \DecisionSet^\CurrentTime(\StateNote_\CurrentTime)} \{ \Reward(\DecisionTuple^\CurrentTime) + V_\CurrentTime^\PostDecision(\StateNote_\CurrentTime^\VehiclePost) \} \label{eq: BellOne} \\
    & V_\CurrentTime^\PostDecision(\StateNote_\CurrentTime^\VehiclePost)=\mathbb{E}_{\ExogenousInformation_\CurrentTime}\big[V_{\CurrentTime + 1}(\StateNote_{\CurrentTime + 1})\big|\StateNote_\CurrentTime^\VehiclePost\big]  \label{eq: BellTwo}
    \end{alignat}
\end{subequations}
Equation~\eqref{eq: BellOne} establishes a deterministic optimality equation based on the post-decision state, eliminating the need to enumerate the complete outcome space and assess future values, while Equation~\eqref{eq: BellTwo} expresses the post-decision state value function as the expected value of downstream rewards, such that $\StateNote_{\CurrentTime + 1}=\SecondTransition(\StateNote_\CurrentTime^\VehiclePost,\ExogenousInformation_\CurrentTime)$. In order to solve the Bellman optimality equation by utilizing post-decision states, it is necessary to have knowledge of the value functions for every feasible post-decision state. However, given the high dimensionality of the post-decision state space, traditional approaches are typically incapable of computing the value functions for every post-decision state. To overcome this challenge, we adopt an approximation of the post-decision state value function, which is provided in the following section.

\subsection{Value Function Approximation}
The post-decision state in the RMP is influenced by the vehicles and their potential states, which are determined by factors such as their locations, the number of assigned passenger groups, the pick-up and drop-off locations for each group, the number of passengers per group, and the decision epoch. This increase in the complexity is referred to as the ``second curse of dimensionality'', which pertains to the dimensionality of the state space. To tackle this issue, we propose an approximation of the post-decision state value function, denoted as $V^{\PostDecision}_\CurrentTime (\StateNote_\CurrentTime^\VehiclePost)$, using a linear decomposition function $\bar{V}^\PostDecision_\CurrentTime (\StateNote_\CurrentTime^\VehiclePost)$ of the vehicle vector attributes present in the post-decision state, i.e., $\{{\StateNote_{\CurrentTime \IndividualVehicle}^\VehiclePost}\}_{\IndividualVehicle \in \VehiclesSet_\CurrentTime}$. The formal definition of the linear approximation of the post-decision value function is
\begin{equation}
    \bar{V}_\CurrentTime^\PostDecision (\StateNote_\CurrentTime^\VehiclePost) = \sum_{\IndividualVehicle \in \VehiclesSet_\CurrentTime} \EstimatedV_{\CurrentTime \IndividualVehicle}^\PostDecision \StateNote_{\CurrentTime \IndividualVehicle}^\VehiclePost \label{eq: LinApprox}
\end{equation}
where $\EstimatedV_{\CurrentTime \IndividualVehicle}^\PostDecision$ is the expected down-stream reward associated with a vehicle being in the post-decision state of $\IndividualVehicle$ at time $\CurrentTime$. This approach presents a considerable computational benefit. Rather than computing value functions for every possible post-decision state at the current time $\CurrentTime$, we need only $\big|\VehiclesSet_\CurrentTime\big|$ variables, which are denoted as $\EstimatedV_{\CurrentTime \IndividualVehicle}^\PostDecision$. In the following section, we will discuss the sampling approach used to estimate the $\EstimatedV_{\CurrentTime \IndividualVehicle}^\PostDecision$ values.
\subsection{ADP Policy}
The final curse of dimensionality is linked to the action space, provided that the number of actions a vehicle may take increases exponentially with an increase of the action space. This can make it difficult to find optimal solutions, given the computational complexity of solving such a problem. We address this by introducing the ADP policy, \ADP. The decision state value function can be written as the function of matching decisions $\DecisionTuple^\CurrentTime$ as follows:
\begin{equation}
    \bar{V}_\CurrentTime^\PostDecision (\StateNote_\CurrentTime^\VehiclePost) = \sum_{\IndividualVehicle \in \VehiclesSet_\CurrentTime, \MovementPoint \in \MovementSet_\CurrentTime} \EstimatedV_{\CurrentTime \IndividualVehicle'}^\PostDecision \MatchingVariable_{\CurrentTime \IndividualVehicle \MovementPoint}
    \label{eq: DecApprox}
\end{equation}
where $\EstimatedV_{\CurrentTime \IndividualVehicle'}^\PostDecision$ represents the value of post-decision vehicle state $\IndividualVehicle'$ which is derived by assigning the vehicle of type $\IndividualVehicle$ the action of type $\MovementPoint$ at time $\CurrentTime$. By subsequently replacing the value function of the post-decision state in \eqref{eq: BellOne} with the linear
approximation introduced in \eqref{eq: DecApprox}, \ADP~decisions are obtained in the following manner:
\begin{equation}
    \DecisionSet_\CurrentTime^\ADP(\StateNote_\CurrentTime) \in \argmax_{\DecisionTuple^\CurrentTime \in \DecisionSet^\CurrentTime(\StateNote_\CurrentTime)} \big\{ \Reward(\DecisionTuple^\CurrentTime) +  \bar{V}_\CurrentTime^\PostDecision (\StateNote_\CurrentTime^\VehiclePost) \big\}
\end{equation}
The subsequent \ADP~policy seeks to match vehicles with passenger requests based upon an evaluation of the marginal reward gained from serving those requests compared to that of being assigned to a different request or to be re-positioned. This matching process is deemed feasible due to the unimodular characteristics of the flow conservation constraints, allowing for an integer solution through the solution of a linear program. Modern optimization solvers are equipped to handle this efficiently, even in problem settings involving significant scale, such as the RMP. The implementation of this policy requires only knowledge of the value function approximations, rather the the value functions themselves. The following section presents the methodology employed to make the aforementioned approximations.

\subsection{Forward-Pass ADP Algorithm}
We next present the forward-pass ADP framework which is used for value function estimation. The original forward-pass algorithm was introduced by \citet{simao2009approximate} for the large-scale fleet management problem and has since been re-purposed in recent RMP studies by \citet{yu2019integrated} and \citet{al2020approximate} for their respective problem settings. The provided derivations, approximations, and algorithms in this section are mostly based on the previous studies, but have been tailored to suit our unique problem setting. To keep track of the algorithmic procedure's main iterations, which corresponds to a new sample path of exogenous information, we index them using the superscript $\SingleDay$ and apply this index to all MDP model and ADP algorithm components. Algorithm~\ref{alg: ADPForward} demonstrates the main steps of the forward-pass algorithm.

\renewcommand\algorithm[1][!ht]{\Algorithm[#1]\setstretch{1.5}}
\begin{algorithm}
\caption{Forward Pass ADP Algorithm}\label{alg:cap} \label{alg: ADPForward}
\begin{algorithmic}
\State \textbf{Step 0:} Initialize post-decision values $(\EstimatedV_{\CurrentTime \IndividualVehicle}^{\PostDecision,0})_{\CurrentTime \in \Time, \IndividualVehicle \in \VehiclesSet_\CurrentTime}$ and step-size $\StepSize^1$
\State \textbf{Step 1:} For $\SingleDay=1, \ldots, \Days$:
\State \hspace{1cm} \textbf{Step 2:} Initialize vehicle locations and choose a sample path $(\hat{\ExogenousInformation}_\CurrentTime^\SingleDay)_{\CurrentTime \in \Time}$
\State \hspace{1cm} \textbf{Step 3:} For $\CurrentTime=1,\ldots,T$:
\State \hspace{2cm} \textbf{Step (3a):} Solve the ADP optimality Equation~\eqref{eq: ADPOpt} to obtain optimal decisions $\DecisionTuple^{\CurrentTime, \SingleDay}$
\State \hspace{2cm} and compute marginal values of Constraints~\eqref{eq: VehicleFlow}, $\MarginalValue_\CurrentTime^\SingleDay=(\MarginalValue_{\CurrentTime \IndividualVehicle}^\SingleDay)_{\forall \IndividualVehicle \in \VehiclesSet_\CurrentTime:\StateNote_{\CurrentTime \IndividualVehicle}^{\PostDecision,\SingleDay} > 0}$
\State \hspace{2cm} \textbf{Step (3b):} Compute $(\MarginalValue_{{\CurrentTime-1}, \IndividualVehicle}^{\PostDecision, {\SingleDay}})_{\IndividualVehicle \in \VehiclesSet_\CurrentTime}$ using \eqref{eq: MargEquation}, and update post-decision value 
\State \hspace{2cm} function estimates of the previous epoch, $(\EstimatedV_{{\CurrentTime-1}, \IndividualVehicle}^{\PostDecision, {\SingleDay-1}})$, via \eqref{eq: MargEquation2}
\State \hspace{2cm} \textbf{Step (3c):} Based on optimal decisions $\DecisionTuple^{\CurrentTime,\SingleDay}$ and random information realizations  
\State \hspace{2cm}$\hat{\ExogenousInformation}_{\CurrentTime}^\SingleDay$, move from the current state $\StateNote_\CurrentTime^\SingleDay$ to the next state $\StateNote_{\CurrentTime+1}^\SingleDay$ using \eqref{eq: RequestTransition}-\eqref{eq: VehicleTransition}
\State \textbf{Step 4:} Return the final approximation values $\EstimatedV^{\PostDecision,\Days}=(\EstimatedV_{\CurrentTime \IndividualVehicle}^{\PostDecision,\Days})_{\CurrentTime \in \Time, \IndividualVehicle \in \VehiclesSet_\CurrentTime}$ for real-time decision making
\end{algorithmic}
\end{algorithm}

The algorithm begins by initializing the value functions of the post-decision states of the vehicles, along with the step-size. Next, for each iteration, $\SingleDay$, we initialize our vehicles to a randomized location and realize a sample path of exogenous information for each decision epoch. This path of exogenous information, denoted by $(\hat{\ExogenousInformation}_\CurrentTime^\SingleDay)_{\CurrentTime \in \Time}$, is representative of the arrival of passenger requests, namely $(\hat{\ExogenousInformation}_\CurrentTime^{\Request,\SingleDay})_{\CurrentTime \in \Time}$. In Step 3, for each decision epoch, we obtain an optimal matching between vehicles and requests, compute the marginal value of making such decisions, and update our value function approximations prior to moving forward in time to the next sate. More specifically, in Step (3a) we begin by obtaining the optimal decisions $\DecisionTuple^{\CurrentTime,\SingleDay}$ of our LP by solving our ADP optimality equation
\begin{equation}
\OptimalityFunction_\CurrentTime(\StateNote_\CurrentTime^\SingleDay):= \max_{\DecisionTuple^\CurrentTime \in \DecisionSet^\CurrentTime(\StateNote_\CurrentTime^\SingleDay)} \big\{\Reward(\DecisionTuple^\CurrentTime) + \bar{V}_\CurrentTime^{\PostDecision,{n-1}} (\StateNote_\CurrentTime^{\VehiclePost,\SingleDay}) \big\}
    \label{eq: ADPOpt}
\end{equation}
based upon the estimated value function of the vehicle post-decision state from iteration $\SingleDay - 1$:
\begin{equation}
    \bar{V}_\CurrentTime^{\PostDecision,{n-1}} (\StateNote_\CurrentTime^{\VehiclePost,\SingleDay}) = \sum_{\IndividualVehicle \in \VehiclesSet} \bar{v}_{\CurrentTime \IndividualVehicle}^{\PostDecision, {\SingleDay -1}} \StateNote_{\CurrentTime \IndividualVehicle}^{\VehiclePost,\SingleDay}
    \label{eq: ADPOpt2}
\end{equation}
After making the matching decisions, we may then calculate the marginal values associated with each visited vehicle attribute, $\MarginalValue_{\CurrentTime}^\SingleDay=(\MarginalValue_{\CurrentTime \IndividualVehicle}^\SingleDay)_{\forall \IndividualVehicle \in \VehiclesSet}$, by finding the partial derivatives of $\MarginalValue_\CurrentTime$ with respect to the right-hand side of the vehicle flow conservation constraints~\eqref{eq: VehicleFlow}: 
\begin{equation}
    \MarginalValue_{\CurrentTime \IndividualVehicle}^\SingleDay:= \frac{\partial \OptimalityFunction_\CurrentTime(\StateNote_\CurrentTime^\SingleDay)}{\partial \StateNote_{\CurrentTime \IndividualVehicle}^{\Vehicle,\SingleDay}} \quad \quad \quad \forall \IndividualVehicle \in \VehiclesSet
    \label{eq: MargEquation}
\end{equation}
where $\MarginalValue_{\CurrentTime \IndividualVehicle}^\SingleDay$ captures the reward associated with having an additional vehicle with an attribute of $\IndividualVehicle$ at decision epoch $\CurrentTime$.
In Step (3b), we first define the equivalent post-decision states associated with the marginal values calculated in Equation~\eqref{eq: MargEquation}. More specifically, we determine the post-decision states from the previous time-step, $\CurrentTime -1$, which correspond to the dual values associated with the vehicle states in the current time. For computing such marginal values, we may use the following:
\begin{equation}
    \MarginalValue_{{\CurrentTime-1}, \IndividualVehicle}^{\PostDecision,\SingleDay}:= \sum_{\IndividualVehicle' \in \VehiclesSet}
    \frac{\partial \OptimalityFunction_\CurrentTime(\StateNote_\CurrentTime^\SingleDay)}{\partial \StateNote_{\CurrentTime, \IndividualVehicle'}^{\Vehicle,\SingleDay}}
    \frac{\partial \StateNote_{\CurrentTime, \IndividualVehicle'}^{\Vehicle,\SingleDay}}{\partial \StateNote_{{\CurrentTime-1}, \IndividualVehicle}^{\VehiclePost,\SingleDay}}\biggl|_{\ExogenousInformation_{\CurrentTime-1}^{\Request}=\hat{\ExogenousInformation}_{\CurrentTime-1}^{\Request,\SingleDay}}
    \label{eq: MargEquation2}
\end{equation}
such that $\frac{\partial \OptimalityFunction_\CurrentTime(\StateNote_\CurrentTime^\SingleDay)}{\partial \StateNote_{\CurrentTime, \IndividualVehicle'}^{\Vehicle,\SingleDay}}=\MarginalValue_{\CurrentTime, \IndividualVehicle'}^\SingleDay$ and $\frac{\partial \StateNote_{\CurrentTime, \IndividualVehicle'}^{\Vehicle,\SingleDay}}{\partial \StateNote_{{\CurrentTime-1}, \IndividualVehicle}^{\VehiclePost,\SingleDay}}$ represents the relative change in the number of vehicles with attribute $\IndividualVehicle'$ in $\StateNote_{\CurrentTime, \IndividualVehicle'}^{\Vehicle,\SingleDay}$ with respect to the number of requests with attribute $\IndividualVehicle$ in ${ \StateNote_{{\CurrentTime-1}, \IndividualVehicle}^{\VehiclePost,\SingleDay}}$. Due to the exogenous information $\ExogenousInformation_{\CurrentTime-1}^{\Request}$, we must account for all realized random transitions in $\hat{\ExogenousInformation}_{\CurrentTime-1}^{\Request}$. This means that if from epoch $\CurrentTime-1$ to $\CurrentTime$, attribute $\IndividualVehicle$ transitions into attribute $\IndividualVehicle'$ due to exogenous information of $\hat{\ExogenousInformation}_{\CurrentTime-1}^{\Request}$, we need only to optimize over the variable $\MarginalValue_{\CurrentTime,\IndividualVehicle'}^{\SingleDay}$, and update $\MarginalValue_{{\CurrentTime-1},\IndividualVehicle}^{\PostDecision,\SingleDay}=\MarginalValue_{\CurrentTime,\IndividualVehicle'}^\SingleDay$. Once the marginal values are otained, we may update the post-decision state value function $\EstimatedV_{{\CurrentTime-1},\IndividualVehicle}^{\PostDecision,{\SingleDay-1}}$ as follows:
\begin{equation}
    \EstimatedV_{{\CurrentTime-1},\IndividualVehicle}^{\PostDecision,{\SingleDay}}=(1-\StepSize^\SingleDay)~\EstimatedV_{{\CurrentTime-1},\IndividualVehicle}^{\PostDecision,{\SingleDay-1}} + \StepSize~\MarginalValue_{{\CurrentTime-1},\IndividualVehicle}^{\PostDecision,\SingleDay}
    \label{eq: UpdateEquation}
\end{equation}
where $\StepSize^\SingleDay$ is the step-size, which we define in further details in the following subsection. Finally, in Step (3c) we transition the state of our system forward in time, and in Step 4, once we have iterated over all of our sample paths, we return the approximated value for each vehicle post-decision state.

\subsection{Algorithmic Enhancements}

We employ two critical algorithmic enhancements in our ADP algorithm, namely, hierarchical aggregation (HA) and bias-adjusted Kalman filter (BAKF) step-size. First introduced by \citet{george2008value}, HA is commonly employed in ADP for problems with large state spaces such as RMP. It combines value function approximations from different levels of aggregation with weighted proportions to improve the approximation of disaggregated values. This technique offers several advantages. Firstly, HA significantly reduces the number of states which need to be considered to reach a viable approximation of value functions, thereby reducing the computational complexity of the problem. Secondly, HA can effectively improve the generalization of the value function approximation by mitigating the influence of noise and outliers in the data. Lastly, the approach can enhance the scalability of the algorithm by restricting the number of states that need to be evaluated, making it possible to solve more complex and larger-scale problems. In our RMP problem setting, we look to perform multi-level aggregation with respect to the spatial properties in our vehicle states by clustering neighbouring intersection locations. Provided that we anticipate neighboring locations to exhibit similar value functions, doing this allows us to gain information about the general values of less frequent nodes in our system. More specifically, in our problem, pickup and drop-off locations situated in regions with high demand activity (or ``hot spots'') are deemed more desirable than those located in less active areas, as they result in a higher number of fulfilled requests. As a consequence, we anticipate that neighboring locations exhibit similar value functions. In practice, this leads to vehicles accepting requests which are located in more active areas, allowing for a greater number of pooling opportunities, as well as rebalancing to these locations, allowing for shorter wait and delay times until requests are fulfilled.

Let $\AggregationSet=\{0,1, \ldots, G\}$ to represent the index set of aggregation levels, such that 0 refers to the most disaggregated level, and $\Locations^\AggregateLevel$ represents the set of locations for the aggregation level $\AggregateLevel \in \AggregationSet$. Additionally, we define $\AggregationFunction^\AggregateLevel:\Locations \rightarrow \Locations^\AggregateLevel$ for $\AggregateLevel \in \AggregationSet$ as the aggregation function which maps the disaggregate attributes to the attribute with aggregation level $\AggregateLevel$. Thus, for each individual location in the disaggregated level, $\Location \in \Locations$, we have that $\Location^\AggregateLevel=\AggregationFunction^\AggregateLevel(\Location)$. Consider an example with 100 intersection locations and 2 levels of aggregation (i.e., $\AggregationSet=\{0,1\}$), such that level $\AggregateLevel=1$ has 5 zones. Then, the cardinality of the set of locations for the two levels will be 100 and 5, respectively, such that the nodes in the least aggregated level, $\AggregateLevel=0$, are clustered into groups of 20 and assigned to a single zone in $\AggregateLevel=1$. Thus, by transferring the aggregated location levels to the vehicle and requests, we obtain aggregated vehicle post-decision states. More specifically, for each level of aggregation denoted by $\AggregateLevel \in \AggregationSet$, we have that $\EstimatedV_{\CurrentTime-1}^{\PostDecision,\AggregateLevel}=(\EstimatedV_{{\CurrentTime-1},\IndividualVehicle^\AggregateLevel}^{\PostDecision,\AggregateLevel})_{\IndividualVehicle^\AggregateLevel \in \VehiclesSet^\AggregateLevel_\CurrentTime}$. The index $\IndividualVehicle^\AggregateLevel \in \VehiclesSet^\AggregateLevel_\CurrentTime$ of each estimated value corresponds to the vehicle state whose current location and assigned pick-up and drop-off points are defined by the aggregate location levels of $\AggregateLevel$. Subsequently, we obtain an enhanced estimate for the post-decision value functions of attribute $\IndividualVehicle$ by computing a weighted sum of estimates across various aggregation levels as
\begin{equation}
    \EstimatedV_{{\CurrentTime-1},\IndividualVehicle}^{\PostDecision}:=\sum_{\AggregateLevel \in \AggregationSet} \Weights_{{\CurrentTime-1},\IndividualVehicle}^{\AggregateLevel}\EstimatedV
    \label{eq: WeightedEstimates}_{{\CurrentTime-1},\IndividualVehicle^\AggregateLevel}^{\PostDecision,\AggregateLevel}
\end{equation}
where $\Weights^{\AggregateLevel}$ represents the weight associated with aggregation level $\AggregateLevel \in \AggregationSet$.

Additionally, we incorporate the BAKF step-size when updating our value function approximations. Introduced by \citet{george2006adaptive}, the BAKF step-size is an adaptive technique that modifies the step-size of the value function approximation update based on the estimated bias. This process helps reduce the impact of bias on the estimates, thereby improving the algorithm's accuracy. Moreover, the BAKF step-size enhances the algorithm's convergence rate by adjusting the step size according to the local curvature of the value function.

Finally, we include additional auxiliary information which describes the overall state of the system in the vehicle post-decision state. This includes information regarding the number of requests which have arrived at the current time-step $\CurrentTime$, as well as the number of vehicles within a proximity of $\MaxWait$ to each agent prior to making any decisions, which we define as the number of ``nearby'' vehicles. Through this information, we can obtain a better understanding of the degree of competitiveness that exists among vehicles for the available passenger requests. In this way, the value of assigning a vehicle to a request increases when requests are relatively scarce, while the value of a request to a vehicle increases when competition among vehicles is higher, allowing us to produce more accurate value function approximations of vehicles with respect to their local environment.

\section{Experimental Setup} \label{experimentalsetup}
The experiments in our study are implemented in Python 3.6.13 and executed on Compute Canada Cedar servers \citep{ComputeCanada}. Version 12.10.0 of IBM ILOG CPLEX Optimization Studio is used to solve the LP models. Below, we present a comprehensive summary of the experimental settings to assess the performance of our \ADP~policy. We begin by detailing the pair of real-world and synthetic datasets used to conduct our experiments, and then proceed to define the benchmark policies used for comparison.
\subsection{Dataset Description}
We consider two datasets in our numerical study, namely, the New York dataset and Chicago dataset, which we describe below.

\subsubsection{New York Dataset}
\begin{figure}[!ht]
\begin{center}
\includegraphics[width=0.29\textwidth]{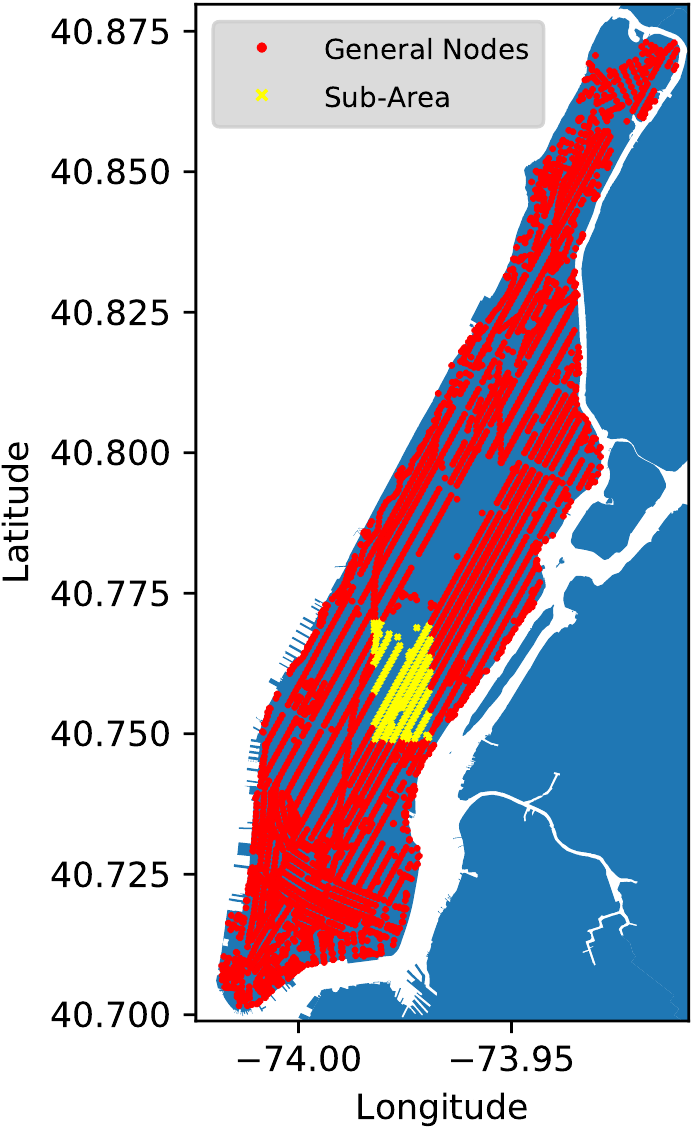}
\end{center}
\caption{Manhattan node locations: nodes in chosen sub-area in yellow (x).}
\label{fig:ny_area}
\end{figure}

\begin{figure}[!ht]
    \centering
    \subfloat[Request distribution\centering]{{\includegraphics[width=0.42\textwidth]{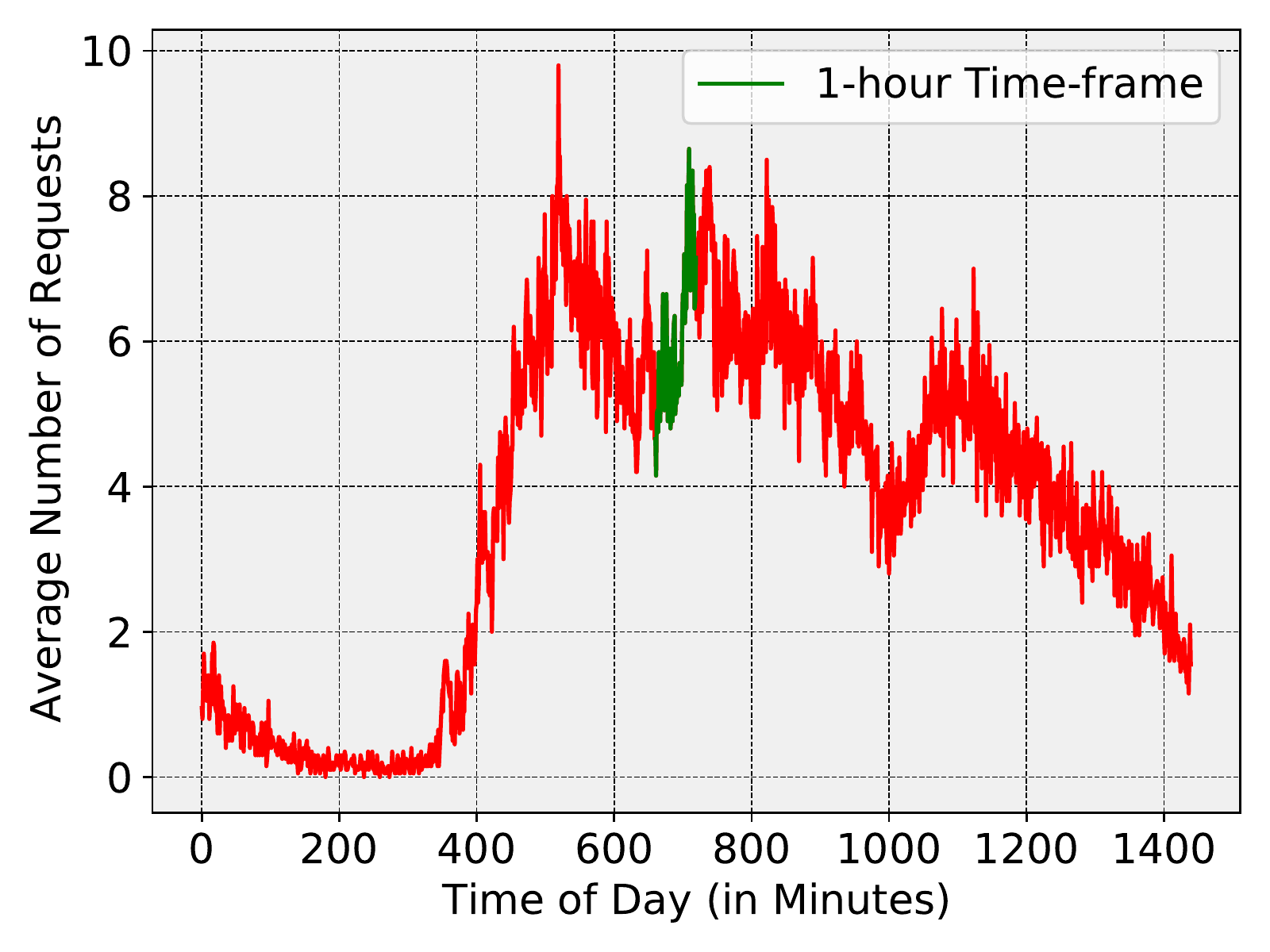} }}
    \subfloat[Passenger size distribution \centering]{{\includegraphics[width=0.42\textwidth]{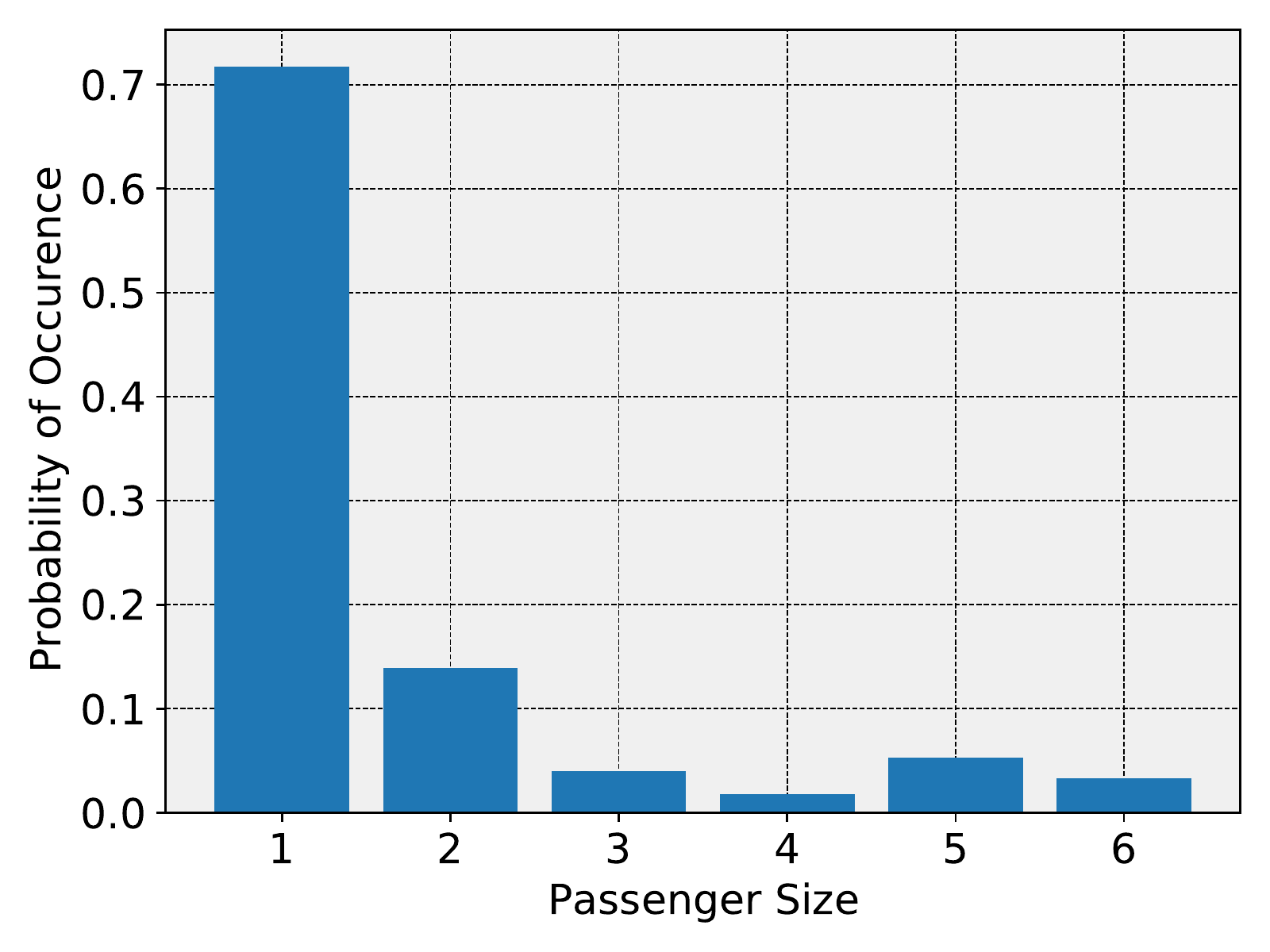} }}%
    \caption{New York dataset specifications.}
    \label{fig:ny_req_dist}
\end{figure}

The \textit{New York Dataset} is a commonly used \citep{alonso2017demand,lowalekar2019zac,shah2020neural,yu2019integrated} comprehensive collection of publicly available \textit{Manhattan Yellow Taxicab} data spanning a period of $20$ days from March to April of 2016~\citep{nyYellowTaxi2016}. It includes information about the daily distribution of ride requests, the popularity of pick-up and drop-off locations, as well as the number of passengers associated with each individual ride. We set this dataset against the backdrop of a real-world Manhattan road network derived using \textit{OpenStreetMap}. In this network, nodes are defined as street intersections and the latitude-longitude coordinates of each pick-up and drop-off point are mapped to the nearest intersection. Edges are defined as the roads connecting each intersection and the weight assigned to each edge represents the travel time between the pair of connected nodes. For the purposes of our experiments, we use a subset of the \textit{Manhattan} area known for high request activity, reducing the number of nodes in the network from 4,442 to 200. The ride request data is then sampled to include only the most frequently visited pairs of pick-up and drop-off locations within this area. A map of this sub-area is provided in Figure~\ref{fig:ny_area}. Furthermore, a time frame of one hour during mid-day rush hour between 11AM and 12PM is selected for the analysis. The overall distribution of requests is derived by taking the average number of requests at each minute of the day over the 20 days worth of data (in the given sub-area), and it is sampled to only include the specified one-hour time frame. The distribution of passenger sizes is derived in a similar manner for the same time frame and area. The resulting distributions are presented in Figure~\ref{fig:ny_req_dist}, in which the specified time frame for the request data is highlighted in green.
\subsubsection{Chicago Dataset}
\begin{figure}[!ht]
\begin{center}
\includegraphics[width=0.7\textwidth]{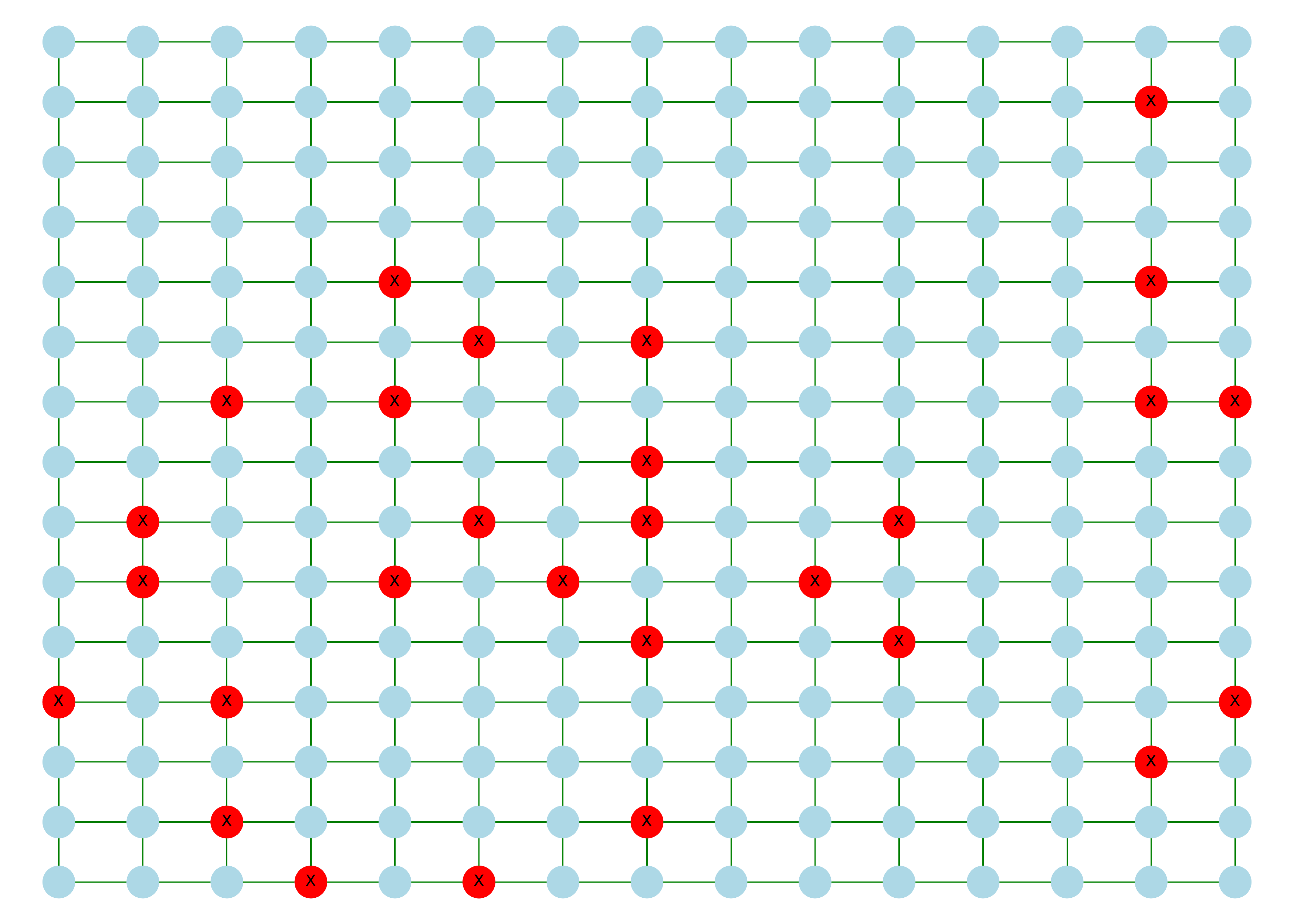}
\end{center}
\caption{Sample Chicago pick-up and drop-off locations in red (x).}
\label{fig:chicago_grid}
\end{figure}

\begin{figure}[!ht]
    \centering
    \subfloat[Request distribution\centering]{{\includegraphics[width=0.42\textwidth]{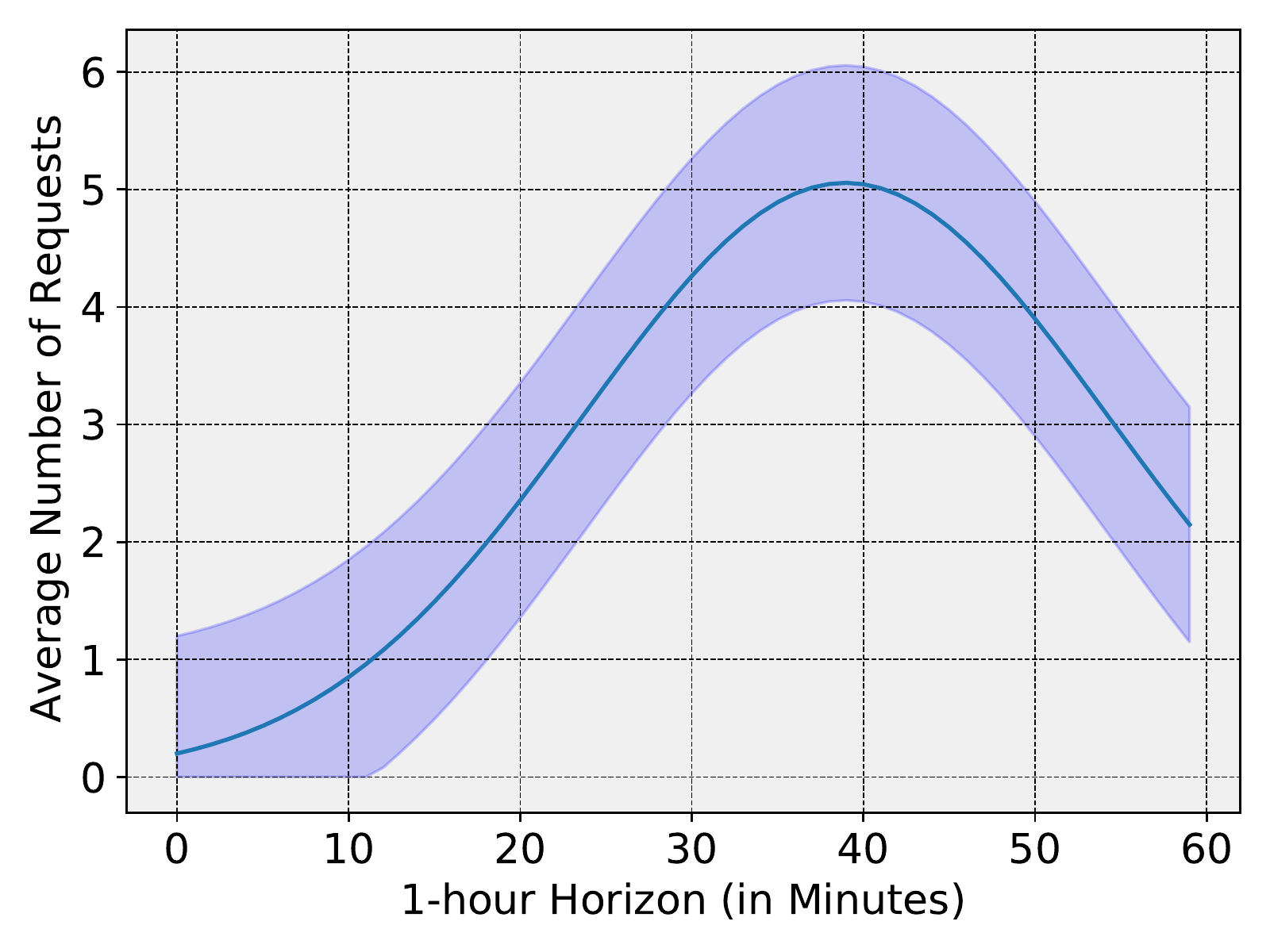} }}%
    \subfloat[Passenger size distribution\centering]{{\includegraphics[width=0.42\textwidth]{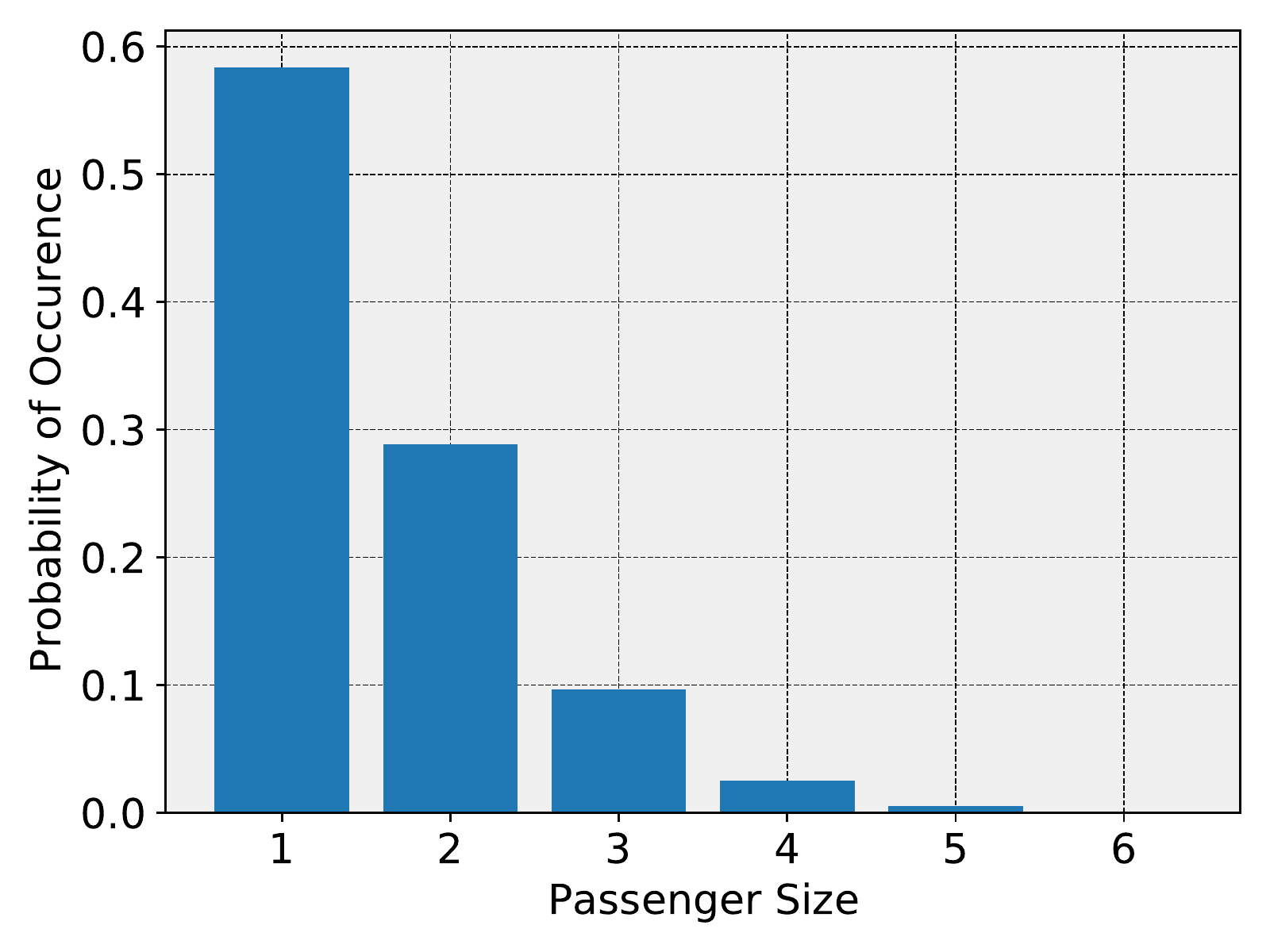} }}%
    \caption{Chicago dataset specifications}
    \label{fig:chicago_req_dist}
\end{figure}

The \textit{Chicago Dataset} is a collection of taxicab requests collected in the city of Chicago between October 1st and 21st of 2019~\citep{chicagoTaxi2019}.
It encompasses information about the daily distribution of ride requests and the prominence of pick-up and drop-off locations. In order to provide a representative sample, the pick-up and drop-off points are selected from the most frequented locations within the dataset. These locations are then clustered using the Louvain algorithm~\citep{Blondel_2008} to preserve community structure and mapped onto a $15\times15$ grid. The grid is broken up into as many zones as there are communities, and the nodes within each community are randomly placed within this sub-zone. In our experiments, both directed and undirected grids are utilized, with arc/edge weights representing travel time ranging between $20$ to $40$ seconds. A map of the grid setup, including an example of the placement of the pick-up and drop-off points, is shown in Figure~\ref{fig:chicago_grid}. The grid layout is selected so as to analyze the algorithmic performance against a more malleable structure, and to test the significance of graph density and complexity. Additionally, the arrival times of requests in the original dataset are partitioned into fifteen-minute intervals. Thus, to accurately simulate the distribution of requests during the 1-hour problem horizon, a truncated normal curve is fit to the request arrival data and taken as the distribution for our experiments with a 1 request standard deviation band. Furthermore, the passenger sizes are also set to follow a truncated normal distribution. A graphical summary of the request and passenger size distributions can be found in Figure~\ref{fig:chicago_req_dist}, in which the curve for the request distribution is shown for the 1-hour horizon.
\subsection{Benchmark Policies}
We consider two benchmark policies to compare our ADP policy against, namely, a Myopic and NeurADP policy, which we describe next.

\subsubsection{Myopic Policy}
Generally, myopic policies are localized or greedy policies which base decisions solely on the current state of the system, disregarding any future consequences of taken actions. This approach simplifies decision-making and allows for quicker execution, making it useful for addressing dynamic and complex issues where computing a globally optimal policy, such as that of the RMP, may not be feasible. However, due to their shortsightedness, myopic policies may not always produce optimal long-term outcomes and are thus often used as baselines. Thus, the myopic benchmark we consider, referred to as \Myopic, prioritizes the number of requests served in the present time step by solving a linear program with constraints~\eqref{eq: VehicleFlow}-\eqref{eq: VariableDefinition}, but does not take into account downstream rewards (i.e., ignores the expected value in the objective function). 
\subsubsection{NeurADP Policy}
We use a Neural Approximate Dynamic Programming policy for our second benchmark policy, which we refer to as \NeurADP. Introduced by \citet{shah2020neural} for a one-to-many case of the RMP, \NeurADP~ is a novel ADP-based algorithm specifically designed to tackle the challenges of large-scale problems where traditional ADP methods fall short. Just as in \ADP~, \NeurADP~ looks to solve sequential decision-making problems by approximating the value functions of post-decision states. However, unlike \ADP~which relies on linear or piece-wise linear value function approximations and utilizes LP-based duals for updating said functions, \NeurADP~ employs a non-linear neural network-based value function approximation and utilizes deep reinforcement learning-based approaches to update them. In particular, the gradients linked with the network parameters are computed and adjusted through the minimization of the L2-norm between the present estimate of the value function and a one-step projection of the return derived from the Bellman equation. Moreover, to enhance stability, off-policy updates are combined with deep reinforcement learning techniques, including the implementation of a target network and Double Q-learning \citep{van2016deep}. At each time step for each vehicle, the current state and feasible action set are stored as experience. During experience replay, the feasible actions are then scored using the value function, and an integer linear program is used to determine the best matching between the saved post-decision state and the generated next post-decision state. The value function of the preserved post-decision state is then updated with the value function of the generated next post-decision state. Regarding the architecture of the neural network, an embedding layer receives input concerning the current location of a vehicle, as well as the assigned pick-up and drop-off points. The embedded representations of these locations, along with their associated delays, are then combined and passed through a long short-term memory (LSTM) layer. Subsequently, this output is concatenated with other relevant auxiliary information and passed through multiple dense layers before producing a single approximate value. Parameter tuning is additionally performed with respect to the overall architecture of the network, including that of the embedding sizes and the number of dense layers.

\section{Results}\label{results}
In this section, we present results from our detailed numerical study. In comparing the performance of the \ADP~policy against the two benchmark policies, we consider the number of requests served within the 1-hour decision horizon. More specifically, we assess each policy by generating five distinct sets of request data spanning 20 days and subsequently compute the average total number of requests that each policy has encountered and fulfilled across these five sets. Passenger requests are generated via sampling at every decision epoch, and a feasible decision is made to match them with vehicles. Upon matching, the vehicle paths are updated, and a simulation of the vehicles' movement along the road network takes place. We first report the findings related to the \textit{New York} dataset, followed by an analysis of the results obtained from the \textit{Chicago} dataset. 
\subsection{Results with the New York Dataset}
The results from the New York data are assessed with respect to four different considerations: wait and delay time, group size, vehicle capacity, and rebalancing. The base case configuration of our parameters are a maximum wait and delay time both set to $90$ seconds, a maximum group size of $3$, and a maximum vehicle capacity of $6$, with the added feature of rebalancing. Additionally, an aggregation size of 3 is used for hierarchical aggregation in the ADP algorithm.
\subsubsection{Impact of Wait and Delay Times}
\begin{figure}[!ht]
\centering
\includegraphics[width=0.60\textwidth]{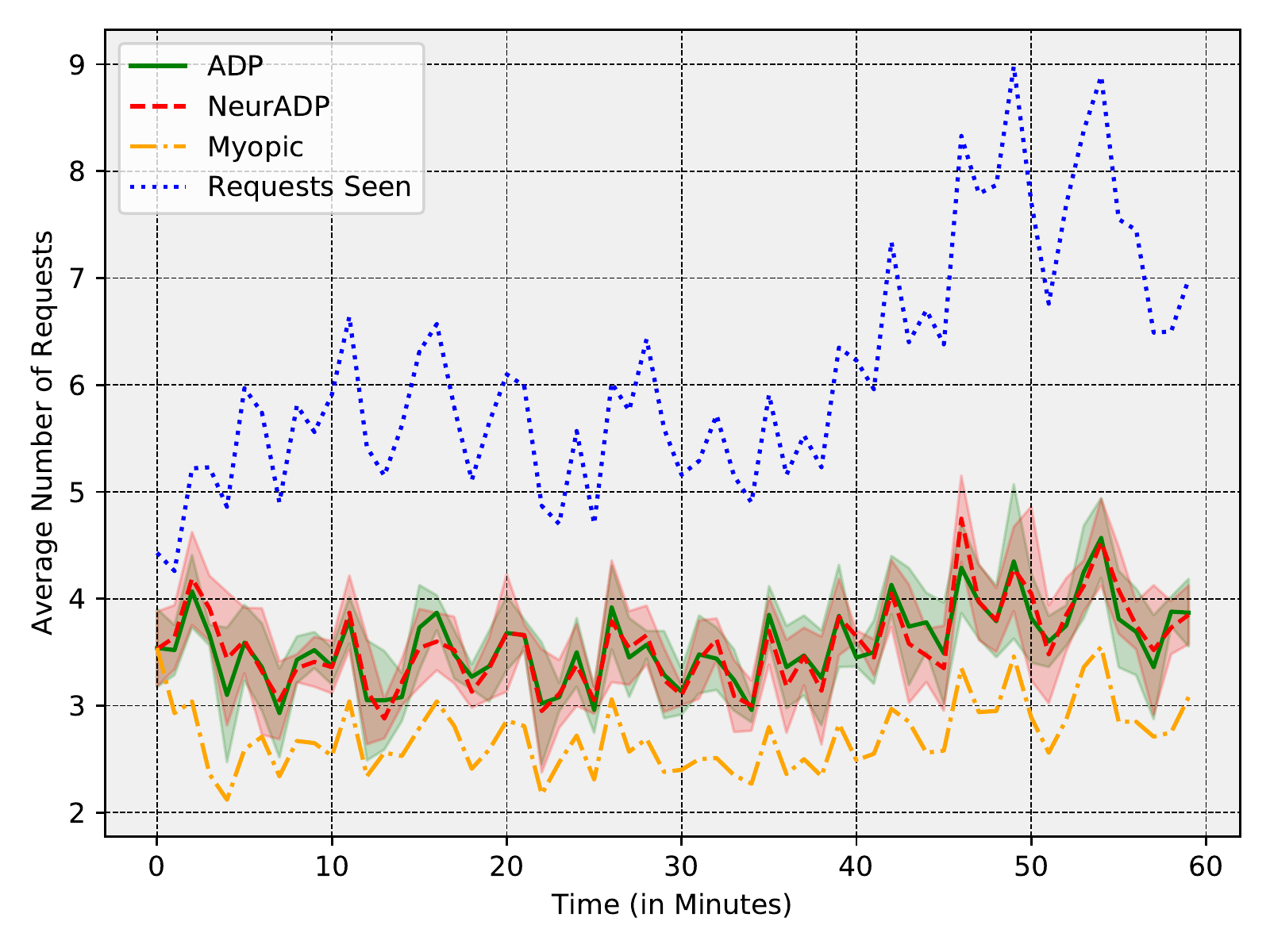}
\caption{Average number of daily requests served (90s W/D time) for the New York dataset.}
\label{fig:ny_daily}
\end{figure}

We first examine how wait and delay times affect the total number of requests served during the operational horizon. Table~\ref{table:delay_table} presents the results for different wait and delay (W/D) times of 60, 90, and 120 seconds, with the wait and delay times being set equal to each other. In this table, the `Requests Seen' column represents the average number of incoming requests observed across the five distinct sets of request data (with `-' denoting a value equal to the row above), while the `\Myopic', `\ADP', and `\NeurADP' columns show the average number of requests served by each policy, with the latter two additionally providing information about the standard deviation of requests served associated with each respective policy. The final two columns illustrate the percentage increase in the average number of requests served by the \ADP~policy compared to the \Myopic~ and \NeurADP~ policies, respectively. This metric is calculated by subtracting the average number of requests served by the baseline policies from the average number of requests served by the \ADP~policy, then dividing the result by the `Requests Seen,' multiplied by 100. These column definitions remain consistent throughout the subsequent tables presented in this paper. Furthermore, Figure~\ref{fig:ny_daily} illustrates the average number of requests served by each policy per decision epoch throughout the 1-hour problem horizon in the scenario with a 90-second delay and wait time. The standard deviation bands are included for each of the non-myopic policies to account for the multiple test days and experiments conducted.

\setlength{\tabcolsep}{4.5pt}
\renewcommand{\arraystretch}{1.15}
\begin{table}[!ht]
\centering
\caption{Impact of Wait \& Delay (W/D) times on the performance as measured by \textbf{average number of requests served} for the New York dataset.}
\label{table:delay_table}
\resizebox{0.95\textwidth}{!}{
\begin{tabular}{P{0.15\textwidth}L{0.12\textwidth}L{0.12\textwidth}C{0.18\textwidth}C{0.18\textwidth}L{0.12\textwidth}L{0.12\textwidth}}
\toprule
\textit{\textbf{W/D time}} & \textit{\textbf{Requests Seen}} & \textit{\textbf{Myopic}} & \textit{\textbf{ADP}} & \textit{\textbf{NeurADP}} & \textit{\textbf{\% Incr. Myopic}} & \textit{\textbf{\% Incr. NeurADP}}\\
\midrule
\textit{\textbf{60 seconds}} & 7,333.4 & 1,365.4 & 3,471.8 $\pm$ 54.1 & 3,469.0 $\pm$ 105.6 & +28.72 & +0.04\\
\textit{\textbf{90 seconds}} & -- & 3,253.2 & 4,280.8 $\pm$ 46.7 &  4,279.6 $\pm$ 61.2 & +14.01 & 0.02\\
\textit{\textbf{120 seconds}} & -- & 4,179.2 & 4,717.8 $\pm$ 20.6 & 4,834.0 $\pm$ 38.4 & +7.34 & -1.58\\
\bottomrule
\end{tabular}
}
\end{table}

\begin{figure}[!ht]
    \centering
    \subfloat[60s W/D time\centering]{{\includegraphics[width=0.42\textwidth]{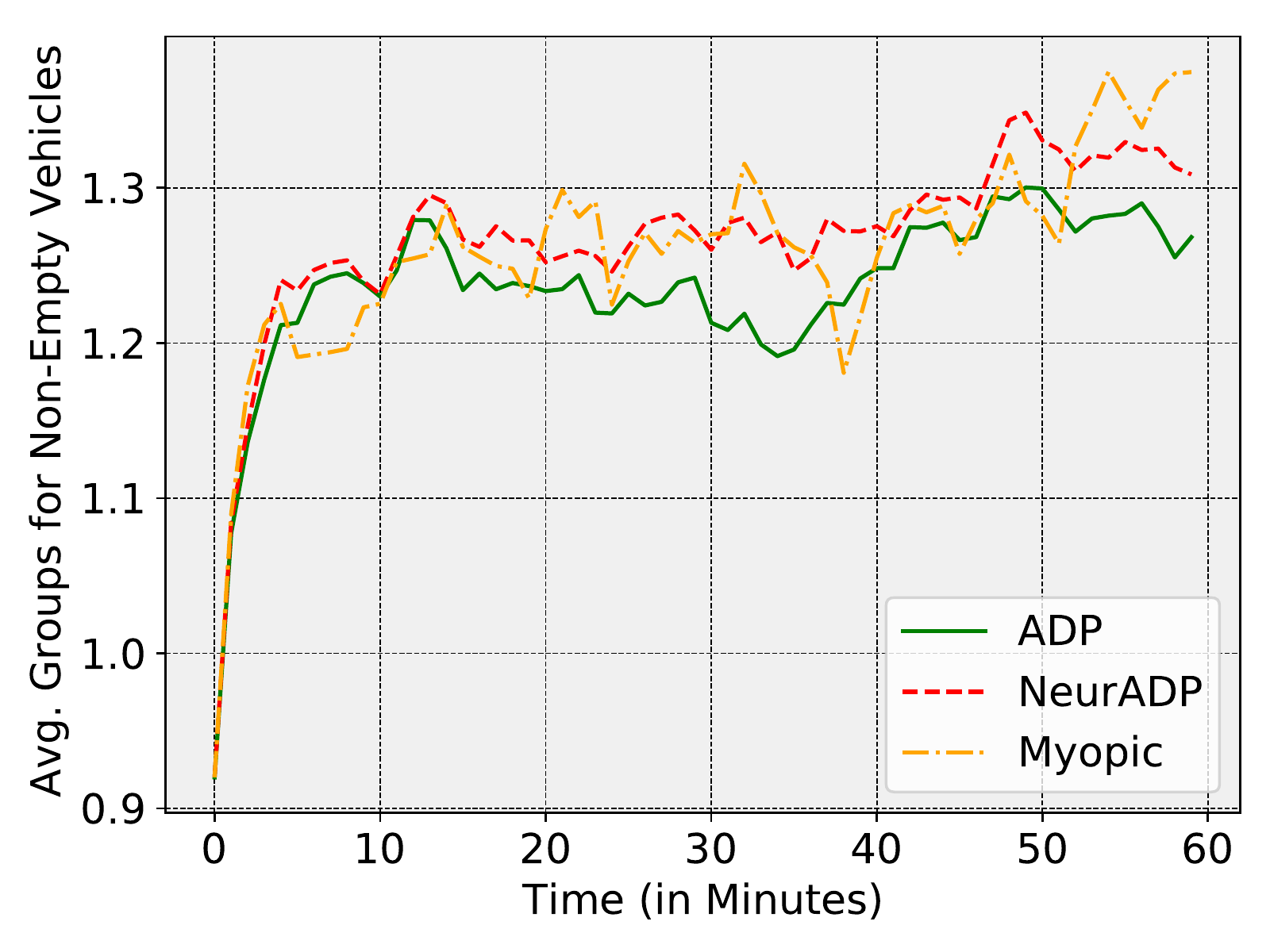}
    }}
    \subfloat[120s W/D time\centering]{{\includegraphics[width=0.42\textwidth]{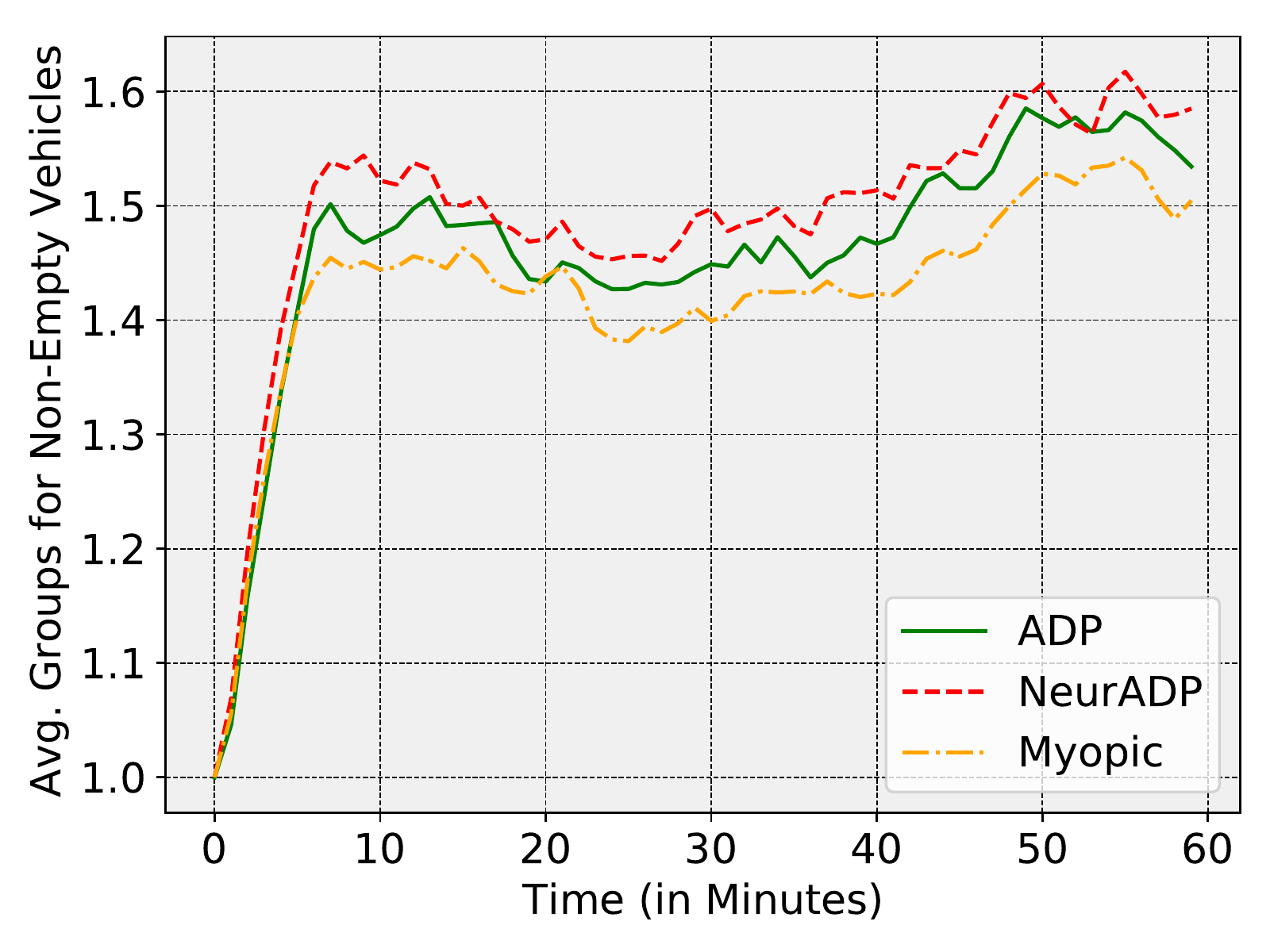}
    }}
    \caption{Passenger group sizes for non-empty vehicles for 60 seconds and 120 seconds wait/delay times for the New York dataset.}
    \label{fig:ny_group_sizes_60_120}
\end{figure}

These results indicate a diminishing rate of increase in requests served in relation to the allowed wait and delay time, irrespective of the employed policy. This is due to the fact that looser time restrictions enable agents to be matched with requests further from their current location and allow for prolonging of existing passengers' rides, thereby enabling the pooling of a higher number of passenger groups. This can be seen in Figure~\ref{fig:ny_group_sizes_60_120}, in which the average number of passenger groups in non-empty vehicles increases with the wait and delay time for the $60$ to $120$ second cases. However, as the pick-up and drop-off time expands, the improvement of the \ADP~over the \Myopic~ policy decreases, reflecting the declining importance of making ``smarter'' decisions as more feasible actions become available. Finally, we see that the \ADP~and \NeurADP~ policies exhibit similar performance across all wait and delay-time variations, with the former demonstrating greater stability and lower standard deviation. This finding suggests that the \ADP~policy may be preferable in practical applications requiring reliable and consistent performance. 

\begin{figure}[!ht]
    \centering
    \subfloat[Number of empty vehicles\centering]{{\includegraphics[width=0.42\textwidth]{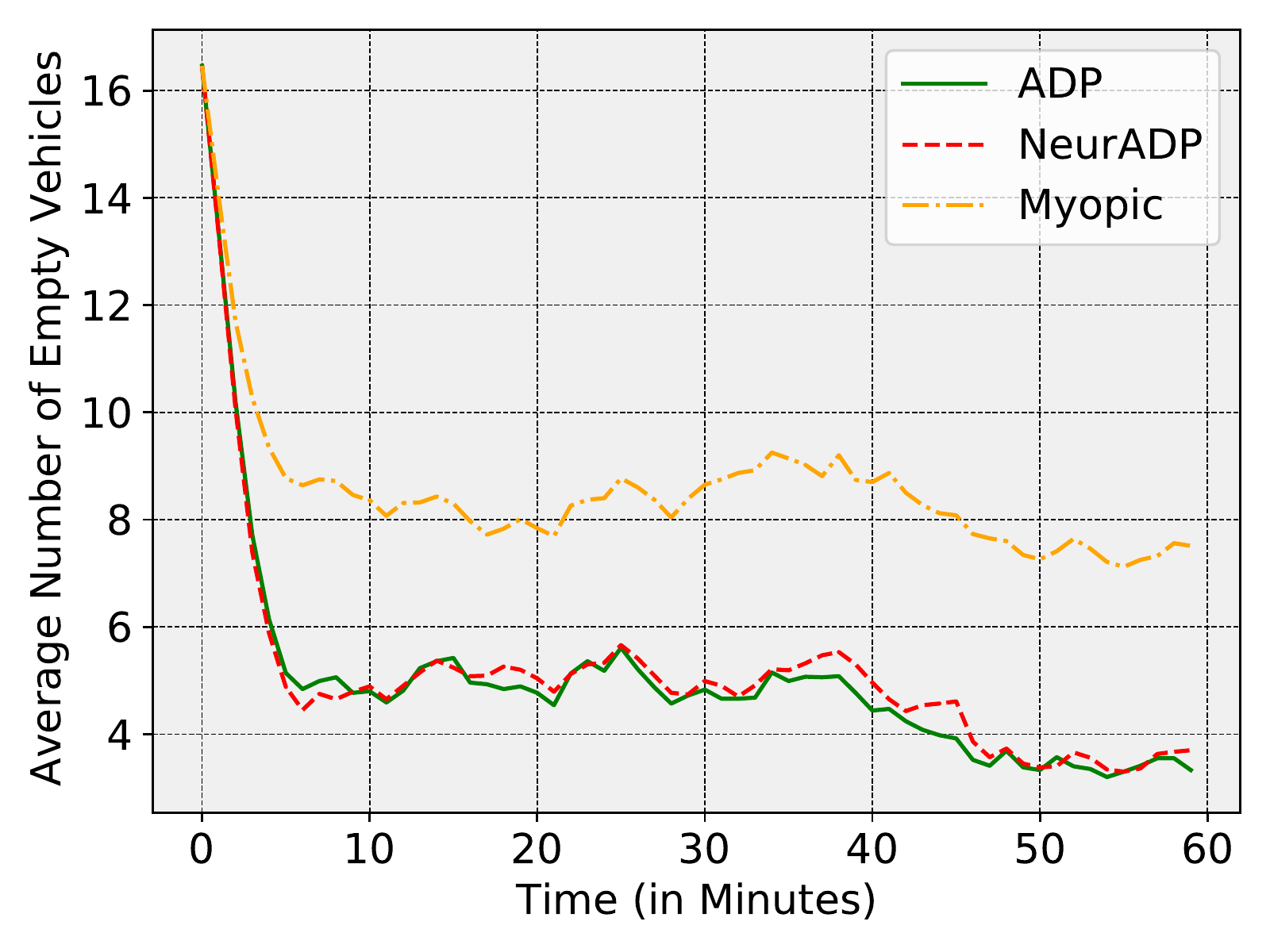}
    }}
    \subfloat[Passenger groups for non-empty vehicles\centering]{{\includegraphics[width=0.42\textwidth]{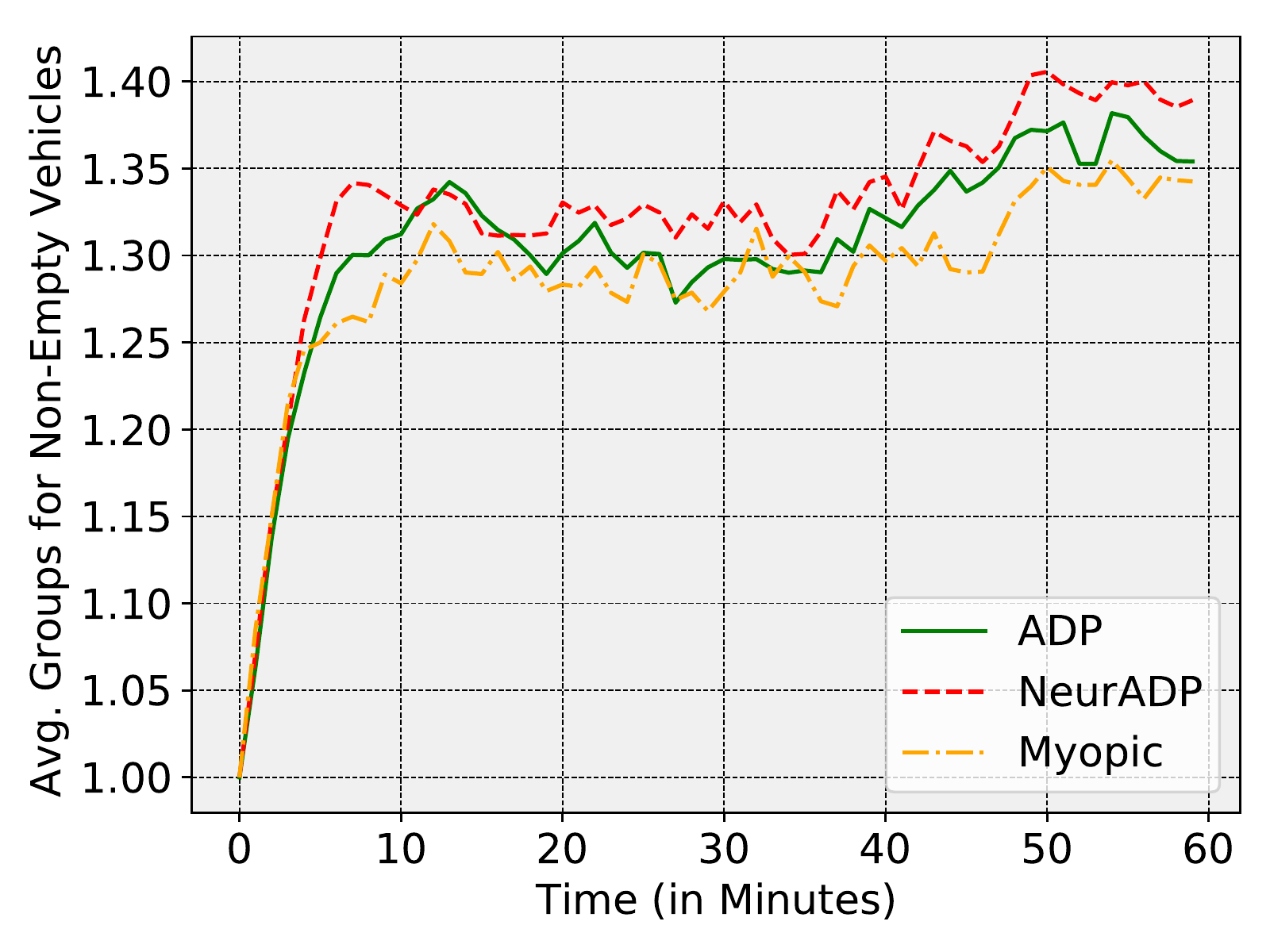}
    }}
    \caption{Vehicle statistics (90s W/D time) for the New York dataset.}
    \label{fig:ny_empty_vehicles}
\end{figure}

We note that, although the \ADP~and \NeurADP~ exhibit comparable performance, they achieve this via differently structured policies. More specifically, the \NeurADP~ policy exhibits a greater tendency to group passengers together, while the \ADP~policy tends to allocate available requests to unoccupied vehicles. This distinction in behavior can be more clearly seen for the 90 seconds case in Figure~\ref{fig:ny_empty_vehicles} in which non-empty vehicles under the \NeurADP~ policy consistently occupy a higher number of passenger groups, while fewer numbers of vehicles under the \ADP~are empty throughout the time-horizon. This contrast may contribute to the observed difference of results in the $120$ second case, as the potential for assigning passenger groups to occupied vehicles is increased.

\subsubsection{Impact of Maximum Group Size}
The relationship between group size and the number of requests served is presented in Table~\ref{table:groups_90_table}. These results indicate a positive correlation between group size and the number of requests served, as larger group sizes allow for the pooling of requests, resulting in more feasible actions and fulfilled requests. However, a diminishing return is observed when comparing the increase in requests served from group size 2 to 3 with that of 1 to 2. This suggests that there may be a limiting factor, such as delay or capacity constraints, that limits the effectiveness of further increasing group sizes on the number of requests served. In order to test this, we conduct a repeated experiment for the 120-second case. Our results indicate that, in comparison to the 90-second case, the increases in performance for the \Myopic, \ADP, and \NeurADP~ requests served are more pronounced when the group size increases from 2 to 3. Specifically, we observe an increase of $1.25\%$, $2.45\%$, and $3.29\%$ for the \Myopic, \ADP, and \NeurADP~ requests served, respectively, for the 120-second case, whereas the corresponding increase is only $0.45\%$, $1.43\%$, and $1.46\%$ for the 90-second case. Furthermore, from Table~\ref{table:groups_90_table} we see a general decrease in improvement from the \ADP~policy over the \Myopic~ policy as the group size increases. This may be caused by the fact that, as the ability to pool a higher number of requests grows, the need for a more sophisticated policy decreases, and as such, the \Myopic~ policy performs more comparatively. Finally, the results indicate that the performance of the \ADP~policy remains comparable to the \NeurADP~ policy across all variations of group size with an improvement in stability.

\setlength{\tabcolsep}{4.5pt}
\renewcommand{\arraystretch}{1.15}
\begin{table}[!ht]
\centering
\caption{Impact of group size (90s W/D time) on the performance as measured by \textbf{average number of requests served} for the New York dataset.}
\label{table:groups_90_table}
\resizebox{0.95\textwidth}{!}{
\begin{tabular}{P{0.16\textwidth}L{0.1\textwidth}L{0.12\textwidth}C{0.18\textwidth}C{0.18\textwidth}L{0.12\textwidth}L{0.12\textwidth}}
\toprule
\textit{\textbf{Group Sizes}} & \textit{\textbf{Requests Seen}} & \textit{\textbf{Myopic}} & \textit{\textbf{ADP}} & \textit{\textbf{NeurADP}} & \textit{\textbf{\% Incr. Myopic}} & \textit{\textbf{\% Incr. NeurADP}}\\
\midrule
\textit{\textbf{1 Groups}} & 7,333.4 & 2,582.4 & 3,371.6 $\pm$ 16.5 & 3,394.0 $\pm$ 22.9 & +10.76 & -0.31\\
\textit{\textbf{2 Groups}} & -- & 3,220.0 & 4,175.6 $\pm$ 30.0 & 4,172.8 $\pm$ 129.8 & +13.03 & +0.04\\
\textit{\textbf{3 Groups}} & -- & 3,253.2 & 4,280.8 $\pm$ 46.7 & 4,279.6 $\pm$ 61.2 & +14.01 & +0.02\\
\bottomrule
\end{tabular}
}
\end{table}

\subsubsection{Impact of Maximum Vehicle Capacity}
\begin{figure}[!ht]
    \centering
    \subfloat[3 passengers\centering]{{\includegraphics[width=0.42\textwidth]{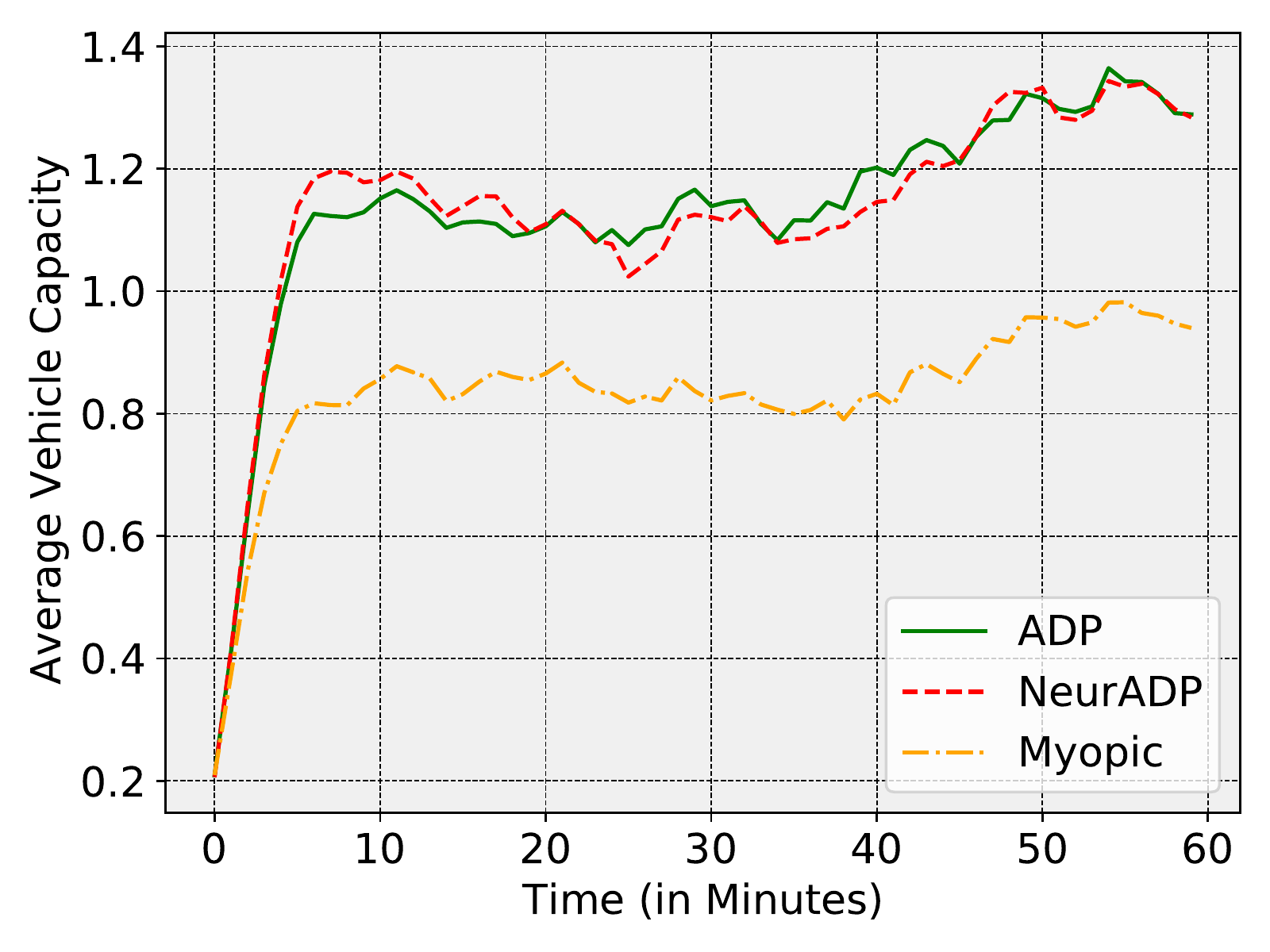}
    }}
    \subfloat[6 passengers\centering]{{\includegraphics[width=0.42\textwidth]{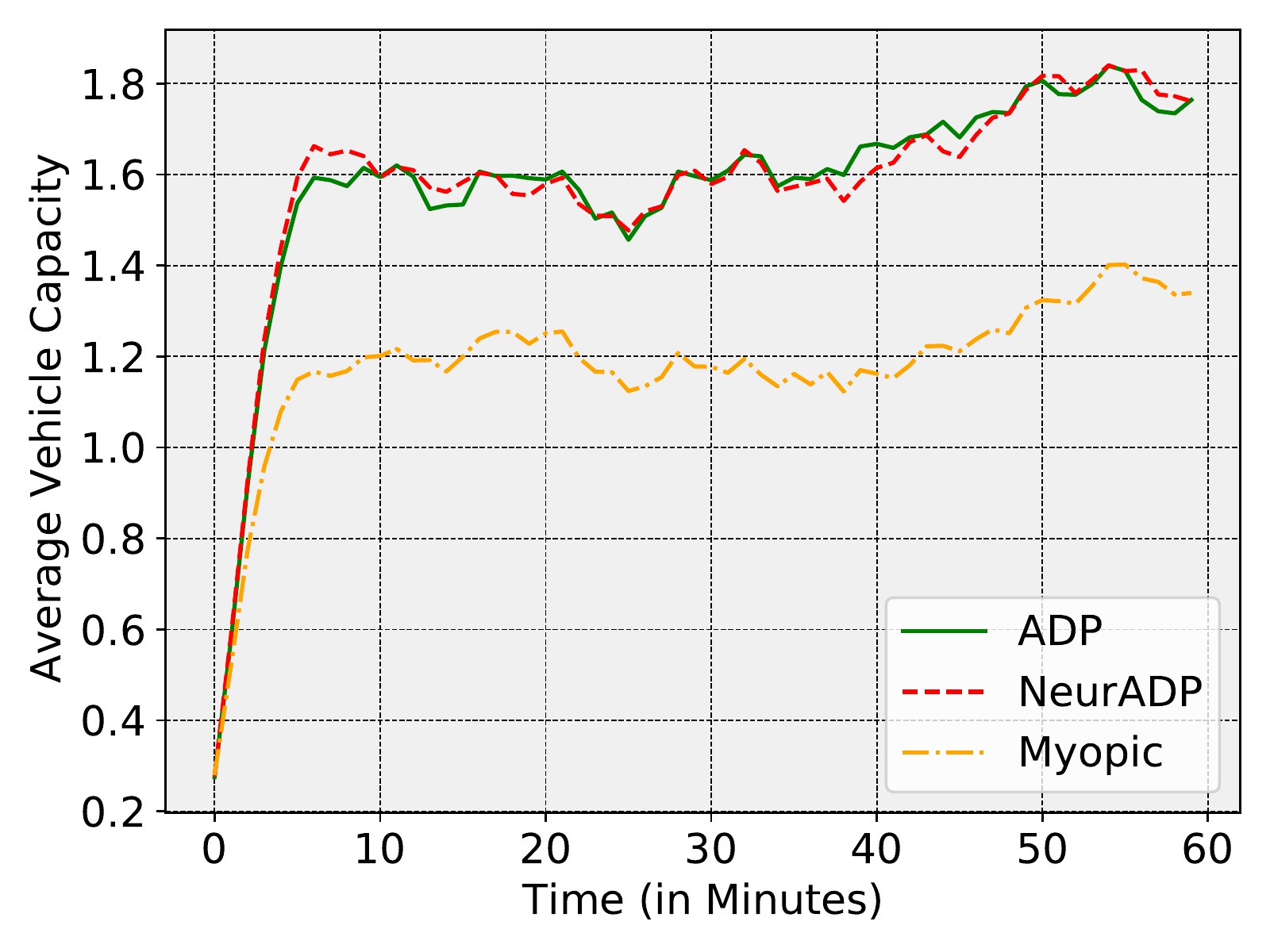}
    }}
    \caption{Average utilized vehicle capacity for different maximum capacity sizes for the New York dataset.}
    \label{fig:ny_vehicle_capacities}
\end{figure}

We next examine the impact of capacity size on the number of requests fulfilled and present the findings in Table~\ref{table:capacity_table}. The results demonstrate that a larger allowed capacity size leads to a proportional increase in the number of requests served due to the availability of more vehicle space to match a greater number of requests. This correlation can be visualized in Figure~\ref{fig:ny_vehicle_capacities}, where the  average utilized vehicle capacity for all policies increases as the maximum allowed capacity increases from 3 to 6. Moreover, we find that the improvement brought by the \ADP~policy over the \Myopic~ policy also grows as the capacity size increases. This is likely due to the ability to accommodate more requests as capacity increases, thus leading to more room for improvement. Finally, our findings indicate that the \ADP~policy remains comparable to the \NeurADP~ policy in terms of performance for different capacity sizes.
\setlength{\tabcolsep}{4.5pt}
\renewcommand{\arraystretch}{1.15}
\begin{table}[!ht]
\centering
\caption{Impact of capacity size (90s W/D time) on the performance as measured by \textbf{average number of requests served} for the New York dataset.}
\label{table:capacity_table}
\resizebox{0.95\textwidth}{!}{
\begin{tabular}{P{0.16\textwidth}L{0.1\textwidth}L{0.12\textwidth}C{0.18\textwidth}C{0.18\textwidth}L{0.12\textwidth}L{0.12\textwidth}}
\toprule
\textit{\textbf{Capacity Sizes}} & \textit{\textbf{Requests Seen}} & \textit{\textbf{Myopic}} & \textit{\textbf{ADP}} & \textit{\textbf{NeurADP}} & \textit{\textbf{\% Incr. Myopic}} & \textit{\textbf{\% Incr. NeurADP}}\\
\midrule
\textit{\textbf{3 Passengers}} & 7,333.4 & 3,012.8 & 3,983.6 $\pm$ 46.9 & 3,978.4 $\pm$ 79.3 & +13.24 & +0.07\\
\textit{\textbf{4 Passengers}} & -- & 3,085.6 & 4,096.4 $\pm$ 54.7 & 4,086.8 $\pm$ 76.5 & +13.78 & +0.13\\
\textit{\textbf{5 Passengers}} & -- & 3,199.6 & 4,192.0 $\pm$ 36.9 & 4,215.4 $\pm$ 70.7 & +13.53 & -0.32\\
\textit{\textbf{6 Passengers}} & -- & 3,253.2 & 4,280.8 $\pm$ 46.7 & 4,279.6 $\pm$ 61.2 & +14.01 & +0.02\\
\bottomrule
\end{tabular}
}
\end{table}

\subsubsection{Rebalancing Results}
\begin{figure}[!ht]
    \centering
    \subfloat[60s W/D time\centering]{{\includegraphics[width=0.42\textwidth]{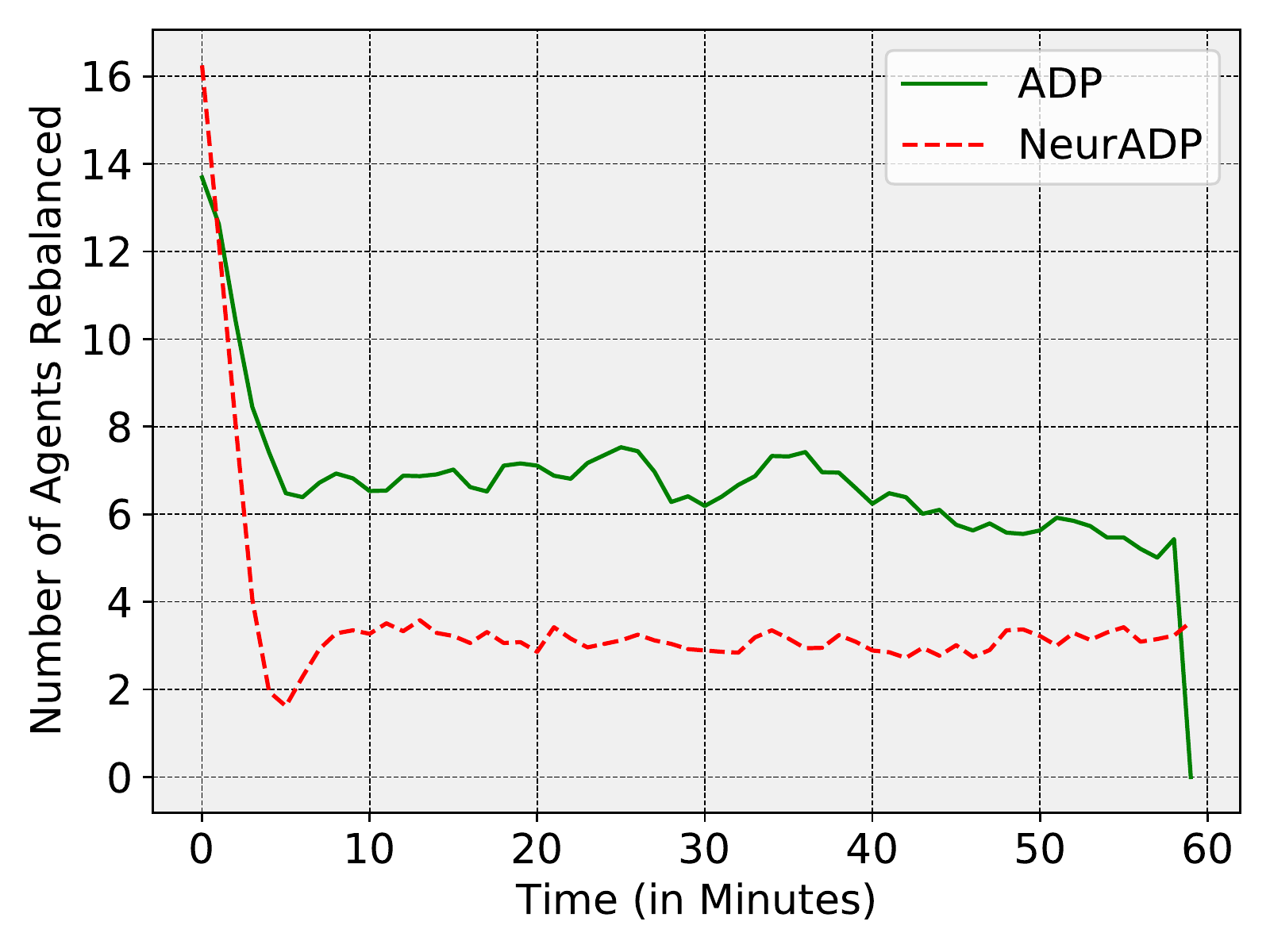}
    }}
    \subfloat[90s W/D time\centering]{{\includegraphics[width=0.42\textwidth]{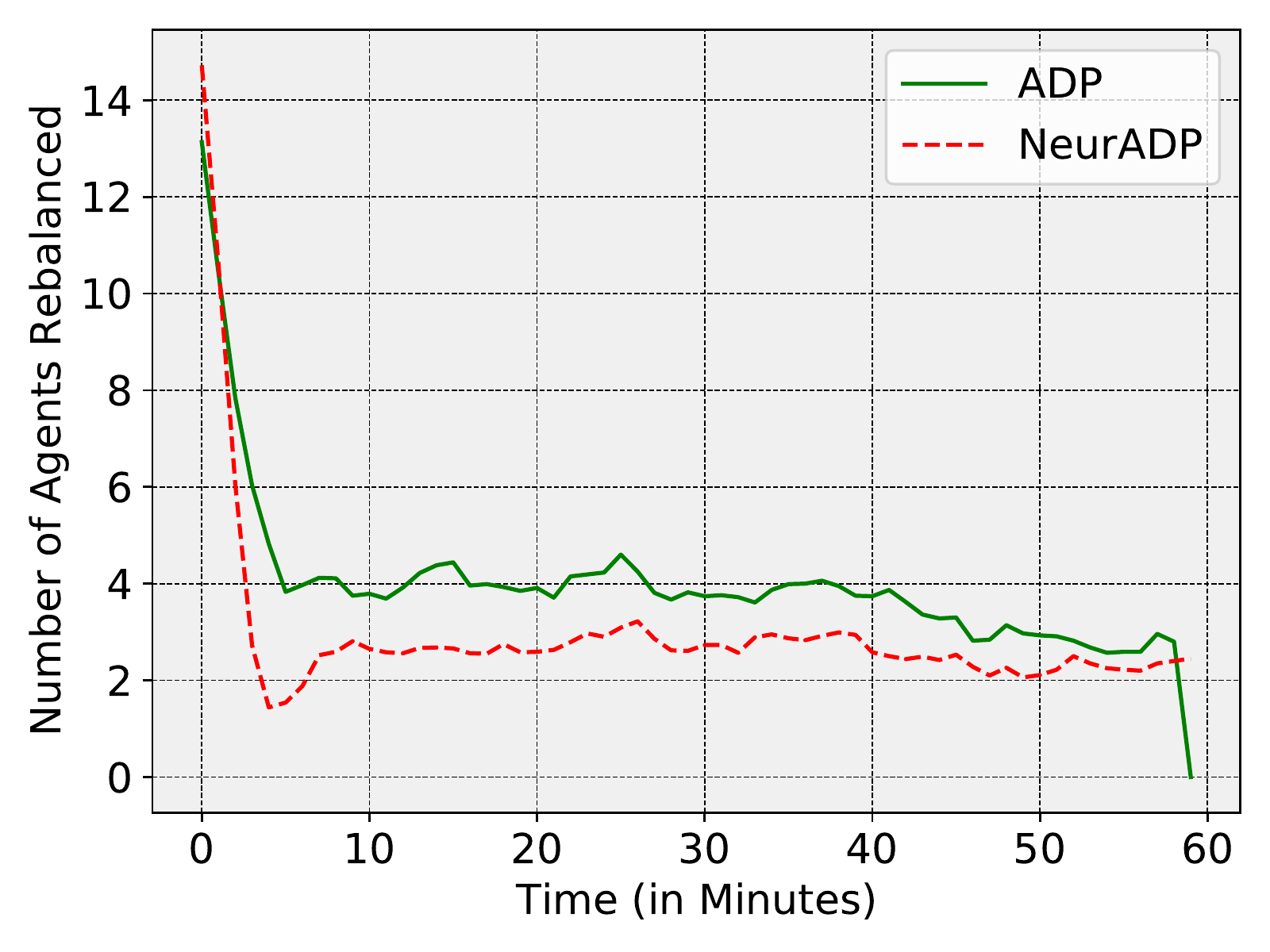}
    }}
    \caption{Average number of vehicles rebalanced for the New York dataset.}
    \label{fig:ny_rebalanced}
\end{figure}

We examine the impact of allowing rebalancing on the number of requests served in Table~\ref{table:rebalancing_table}. Our findings reveal that, in general, allowing rebalancing leads to a higher number of served requests across different policies and wait/delay times. This is due to the fact that rebalancing allows unoccupied vehicles to preemptively travel to popular pickup points, which in turn allows for quicker pick-up times and enables faster trips and more time to pool passengers. Furthermore, our findings indicate that the \NeurADP~ policy surpasses the \ADP~policy in non-rebalancing conditions. This discrepancy from the rebalancing outcomes could be attributed to the \ADP~policy's heightened reliance on rebalancing, as depicted in Figure~\ref{fig:ny_rebalanced}. As a result, the absence of rebalancing has a greater impact on the performance of \ADP~compared to \NeurADP.
\setlength{\tabcolsep}{4.5pt}
\renewcommand{\arraystretch}{1.15}
\begin{table}[!ht]
\centering
\caption{Impact of rebalancing (90s W/D time) on the algorithmic performance for the New York dataset.}
\label{table:rebalancing_table}
\resizebox{0.95\textwidth}{!}{
\begin{tabular}{P{0.25\textwidth}L{0.1\textwidth}L{0.12\textwidth}C{0.18\textwidth}C{0.18\textwidth}L{0.12\textwidth}L{0.12\textwidth}}
\toprule
\textit{\textbf{Rebalancing}} & \textit{\textbf{Requests Seen}} & \textit{\textbf{Myopic}} & \textit{\textbf{ADP}} & \textit{\textbf{NeurADP}} & \textit{\textbf{\% Incr. Myopic}} & \textit{\textbf{\% Incr. NeurADP}}\\
\midrule
\textit{\textbf{Rebalancing (60s)}} & 7,333.4 & 1,365.4 & 3,471.8 $\pm$ 54.1 & 3,469.0 $\pm$ 105.6 & +28.72 & +0.04\\
\textit{\textbf{Rebalancing (90s)}} & -- & 3,253.2 & 4,280.8 $\pm$ 46.7 & 4,279.6 $\pm$ 61.2 & +14.01 & +0.02\\
\textit{\textbf{Rebalancing (120s)}} & -- & 4,179.2 & 4,717.8 $\pm$ 20.6 & 4,834.0 $\pm$ 38.4 & +7.34 & -1.58\\
\midrule
\textit{\textbf{No Rebalancing (60s)}} & -- & 1,365.4 & 1,371.2 $\pm$ 62.5 & 1,524.2 $\pm$ 55.2 & +0.08 & -2.09\\
\textit{\textbf{No Rebalancing (90s)}} & -- & 3,253.2 & 3,326.8 $\pm$ 90.8 & 3,608.0 $\pm$ 117.7 & +1.00 & -3.83\\
\textit{\textbf{No Rebalancing (120s)}} & -- & 4,179.2 & 4,296.2 $\pm$ 67.5 & 4,563.4 $\pm$ 67.9 & +1.60 & -3.64\\
\bottomrule
\end{tabular}
}
\end{table}

\subsection{Results with the Chicago Dataset}
\begin{figure}[!ht]
\centering
\includegraphics[width=0.60\textwidth]{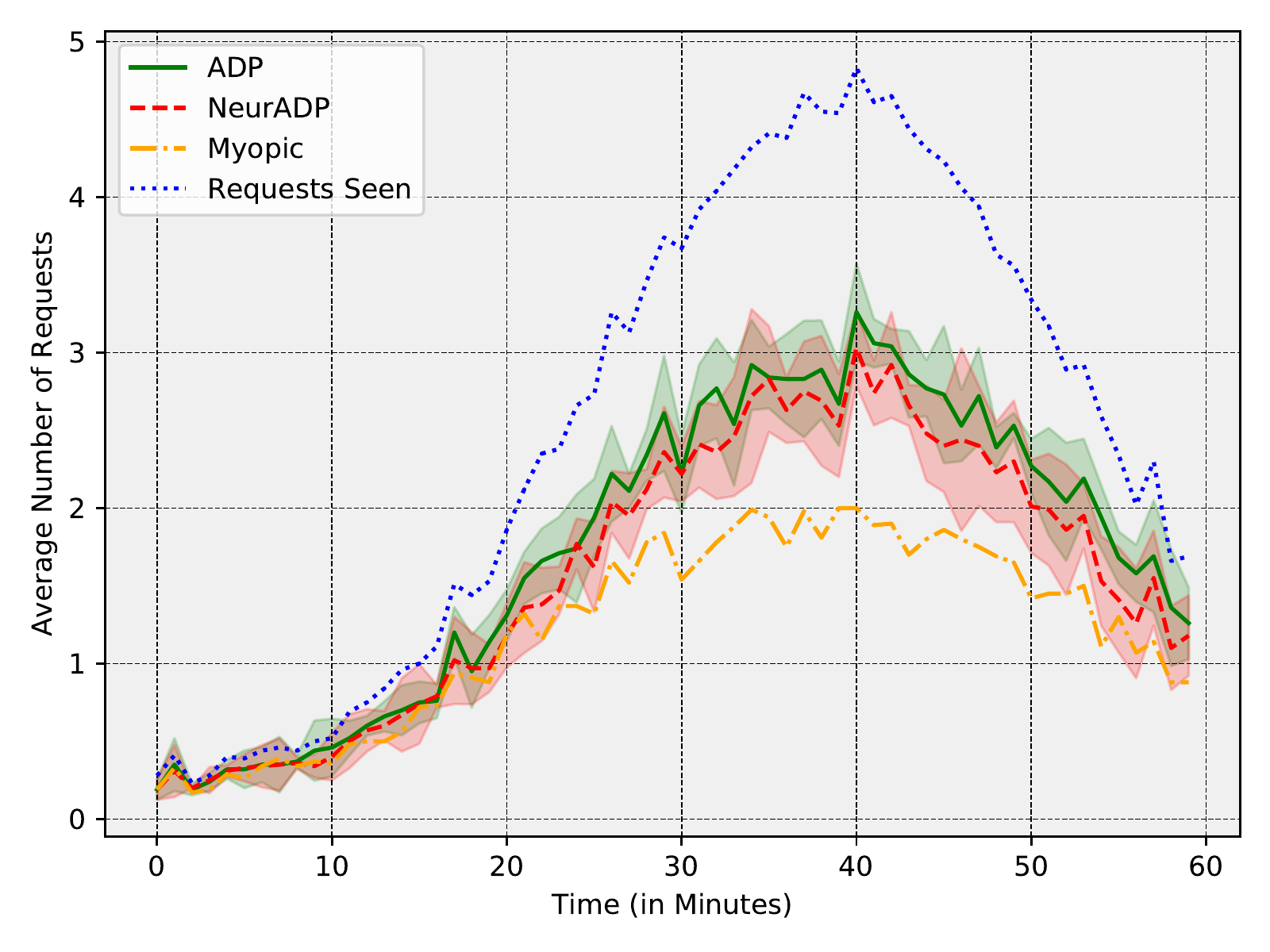}
\caption{Average daily requests served (20\% of edges removed and 100s W/D time) for the Chicago dataset.}
\label{fig:chicago_daily_results}
\end{figure}

The \textit{Chicago} dataset results are divided into two categories: undirected and directed networks. The simulations for the undirected graphs are conducted on road networks that do not have a specified direction, while the simulations for the directed networks are performed on networks with specific road directions. Within each category, the impact of network density on the policy quality is assessed by eliminating portions of the edges (roads) in the graph. This approach not only generates a supplementary dataset for testing the policies but also enables an evaluation of the significance of network complexity. In the undirected network experiments, the edge density of the network is decreased by 20\% and 40\%, while in the directed network experiments, they are reduced by 30\% and 40\%. In both cases, the wait and delay times of 60, 80, and 100 seconds are considered, which are slightly lower than those used for the New York network due to the shorter edge lengths in the grid setup. Table~\ref{table:chicago_undirected_table} and Table~\ref{table:chicago_directed_table} present the results for the undirected and directed network cases, respectively. Additionally, a sample of the daily request service results for the undirected 20\% edge reduction scenario with a wait and delay time of 100 seconds is shown in Figure~\ref{fig:chicago_daily_results}. 

We observe that, in both undirected and directed cases, the number of requests served is again strongly correlated with wait and delay times, with an increase in the pair leading to an increase in the number of requests served. Additionally, the reduction of edges in the network leads to a sharp decline in the number of requests served, as vehicles struggle to traverse the network to meet deadlines and previously feasible requests become infeasible. With regards to directionality, the \ADP~policy consistently outperforms the \NeurADP~ policy in all evaluated undirected scenarios, while we see more comparable performance between the two in the directed case. The superior performance of the \ADP~policy in the undirected scenario may be attributed to a reduction in network complexity, which may facilitate the \ADP~policies' ability to learn and integrate information from the environment and exogenous data. Incorporation of directionality in the graph enhances the intricacy of the network, enabling more competitive performance between the two policies. Finally, we see that the stability of the \ADP~policy is consistently superior or comparable to the \NeurADP~ policy.

\setlength{\tabcolsep}{4.5pt}
\renewcommand{\arraystretch}{1.15}
\begin{table}[!ht]
\centering
\caption{Undirected network results for the Chicago dataset.}
\label{table:chicago_undirected_table}
\resizebox{0.95\textwidth}{!}{
\begin{tabular}{P{0.16\textwidth}P{0.16\textwidth}L{0.1\textwidth}L{0.12\textwidth}C{0.18\textwidth}C{0.18\textwidth}L{0.12\textwidth}L{0.12\textwidth}}
\toprule
\textit{\textbf{W/D times}} & \textit{\textbf{\% Edges Removed}} & \textit{\textbf{Requests Seen}} & \textit{\textbf{Myopic}} & \textit{\textbf{ADP}} & \textit{\textbf{NeurADP}} & \textit{\textbf{\% Incr. Myopic}} & \textit{\textbf{\% Incr. NeurADP}}\\
\midrule
\textit{\textbf{60 seconds}} & \textit{\textbf{0}} & 3,034.6 & 867.8 & 1,773.0 $\pm$ 71.9 & 1,669.2 $\pm$ 132.3 & +29.83 & +3.42\\
\textit{\textbf{60 seconds}} & \textit{\textbf{20}} & -- & 740.4 & 1,554.2 $\pm$ 103.0 & 1,527.2 $\pm$ 125.8 & +26.82 & +0.89\\
\textit{\textbf{60 seconds}} & \textit{\textbf{40}} & -- & 416.8 & 1,197.8 $\pm$ 150.0 & 1,105.2 $\pm$ 165.5 & +25.74 & +3.05\\
\cmidrule(l){2-8}
\textit{\textbf{80 seconds}} & \textit{\textbf{0}} & -- & 1,199.2 & 2,016.6 $\pm$ 72.8 & 1,924.2 $\pm$ 121.9 & +24.94 & +3.04\\
\textit{\textbf{80 seconds}} & \textit{\textbf{20}} & -- & 1,087.6 & 1,795.8 $\pm$ 121.3 & 1,695.8 $\pm$ 143.7 & +23.34 & +3.30\\
\textit{\textbf{80 seconds}} & \textit{\textbf{40}} & -- & 617.0 & 1,351.8 $\pm$ 139.4 & 1,193.8 $\pm$ 191.1 & +24.21 & +5.21\\
\cmidrule(l){2-8}
\textit{\textbf{100 seconds}} & \textit{\textbf{0}} & -- & 1,570.2 & 2,273.4 $\pm$ 94.0 & 2,067.0 $\pm$ 231.0 & +23.17 & +6.80\\
\textit{\textbf{100 seconds}} & \textit{\textbf{20}} & -- & 1,450.8 & 2,065.8 $\pm$ 95.4 & 1,890.0 $\pm$ 167.2 & +20.27 & +5.79\\
\textit{\textbf{100 seconds}} & \textit{\textbf{40}} & -- & 826.4 & 1,520.2 $\pm$ 183.6 & 1,407.6 $\pm$ 162.9 & +22.86 & +3.71\\
\bottomrule
\end{tabular}
}
\end{table}

\setlength{\tabcolsep}{4.5pt}
\renewcommand{\arraystretch}{1.15}
\begin{table}[!ht]
\centering
\caption{Directed network results for the Chicago dataset.}
\label{table:chicago_directed_table}
\resizebox{0.95\textwidth}{!}{
\begin{tabular}{P{0.16\textwidth}P{0.16\textwidth}L{0.1\textwidth}L{0.12\textwidth}C{0.18\textwidth}C{0.18\textwidth}L{0.12\textwidth}L{0.12\textwidth}}
\toprule
\textit{\textbf{W/D times}} & \textit{\textbf{\% Edges Removed}} & \textit{\textbf{Requests Seen}} & \textit{\textbf{Myopic}} & \textit{\textbf{ADP}} & \textit{\textbf{NeurADP}} & \textit{\textbf{\% Incr. Myopic}} & \textit{\textbf{\% Incr. NeurADP}}\\
\midrule
\textit{\textbf{60 seconds}} & \textit{\textbf{0}} & 3,034.6 & 892.8 & 1,744.0 $\pm$ 160.2 & 1,754.4 $\pm$ 179.3 & +28.05 & -0.34\\
\textit{\textbf{60 seconds}} & \textit{\textbf{30}} & -- & 586.4 & 1,395.0 $\pm$ 88.3 & 1,406.4 $\pm$ 86.5 & +26.65 & -0.38\\
\textit{\textbf{60 seconds}} & \textit{\textbf{60}} & -- & 289.2 & 680.4 $\pm$ 190.5 & 828.2 $\pm$ 204.5 & +12.89 & -4.87\\
\cmidrule(l){2-8}
\textit{\textbf{80 seconds}} & \textit{\textbf{0}} & -- & 1,235.6 & 1,983.4 $\pm$ 125.7 & 1,849.0 $\pm$ 152.6 & +24.64 & +4.43\\
\textit{\textbf{80 seconds}} & \textit{\textbf{30}} & -- & 829.2 & 1,576.4 $\pm$ 68.5 & 1,606.8 $\pm$ 111.3 & +24.62 & -1.00\\
\textit{\textbf{80 seconds}} & \textit{\textbf{60}} & -- & 371.0 & 797.2 $\pm$ 141.8 & 878.6 $\pm$ 207.1 & +14.04 & -2.68\\
\cmidrule(l){2-8}
\textit{\textbf{100 seconds}} & \textit{\textbf{0}} & -- & 1,602.4 & 2,175.2 $\pm$ 115.8 & 1,971.8 $\pm$ 195.7 & +18.88 & +6.70\\
\textit{\textbf{100 seconds}} & \textit{\textbf{30}} & -- & 1,134.4 & 1,809.4 $\pm$ 113.7 & 1,689.2 $\pm$ 151.2 & +22.24 & +3.96\\
\textit{\textbf{100 seconds}} & \textit{\textbf{60}} & -- & 464.8 & 844.6 $\pm$ 204.5 & 907.0 $\pm$ 205.0 & +12.52 & -2.06\\
\bottomrule
\end{tabular}
}
\end{table}

\subsection{Summary of Experimental Results}
Overall, our general findings indicate that looser time restrictions lead to an increase in the number of requests served, as agents are able to meet the deadlines and in turn be matched to a larger set of passenger requests. The performance of the \ADP~and \NeurADP~ policies is comparable, with the former exhibiting greater stability and lower standard deviation, while the latter has a greater tendency to group passengers together. The results also indicate a positive correlation between group size and the number of requests served, but a diminishing return is observed when comparing the increase in requests served from group size 2 to 3 with the increase of group size from 1 to 2. Moreover, the results demonstrate that a larger capacity size leads to a proportional increase in the number of requests served. Additionally, the Chicago dataset results demonstrate that a reduction in network density leads to a sharp decrease in the number of requests served. Additionally, the \ADP~policy outperforms the \NeurADP~ policy in all evaluated undirected scenarios, with more competitive performance in the directed case. Our numerical results also suggest that network complexity plays a noticeable role in policy performance, with simplified networks favoring the \ADP~policy.
\section{Conclusion}\label{conclusion}
This paper presents a comprehensive approach to address the challenges of maximizing efficiency in ride-pooling services. By expanding upon the state-of-the-art ADP methodology, introducing key extensions, and performing sensitivity analysis on various parameters, we provide insights into how to optimize the provision of ride-pooling services for the benefit of both passengers and corporations. Our study includes the first comparative analysis between traditional ADP and novel NeurADP approaches, based on experiments conducted on real-world New York taxicab data. Such analysis highlights the strengths and limitations of both methodologies. Moreover, our experiments on a novel Chicago taxi-cab dataset highlight the importance of network density and road directionality in policy performance. NeurADP's utilization of neural networks enables it to learn value function approximations more efficiently, making it capable of handling larger and more complex spatial and temporal problem settings than our ADP approach, which is a limitation of our proposed solution methodology. Future work may focus on combining techniques in both ADP and NeurADP approaches to create a hybrid solution methodology, incorporating the use of hierarchical aggregation with the use of neural network value function approximations. Overall, our findings provide valuable insights for practitioners seeking to maximize the efficiency and effectiveness of ride-pooling services, which can contribute to the reduction of emissions and alleviate traffic congestion.

\section*{Disclosure statement}
No potential conflict of interest was reported by the authors.

\section*{Data Availability Statement}
Data is available upon request.


\singlespacing
\bibliographystyle{elsarticle-harv}

\end{document}